\pdfoutput=1
\RequirePackage{ifpdf}
\ifpdf 
\documentclass[pdftex]{sigma}
\else
\documentclass{sigma}
\fi

\numberwithin{equation}{section}

\newtheorem{Theorem}{Theorem}[section]
\newtheorem*{Theorem*}{Theorem}

\newtheorem{Lemma}[Theorem]{Lemma}
\newtheorem{Proposition}[Theorem]{Proposition}

\theoremstyle{definition}

\usepackage{mathdots}
\usepackage{color}
\usepackage[dvipsnames]{xcolor}

\newcommand{\RR}{\mathbb{R}}
\newcommand{\NN}{\mathbb{N}}
\newcommand{\ficone}{{\widehat{1}}}
\newcommand{\rs}{{\color{red}{\star}}}
\newcommand{\nihat}{{\widehat{n_i}}}
\newcommand{\rel}{\text{rel}}
\newcommand{\tr}{\mathbf{t}}
\newcommand{\nn}{\mathbf{n}}
\newcommand{\vv}{\mathbf{v}}
\newcommand{\ee}{\mathbf{e}}
\newcommand{\two}{\vec{\mathbf{2}}}
\newcommand{\one}{\vec{\mathbf{1}}}
\newcommand{\zero}{\vec{\mathbf{0}}}
\newcommand{\LH}{{\widehat{\Lambda}}}

\begin{document}

\allowdisplaybreaks

\newcommand{\arXivNumber}{2306.16251}

\renewcommand{\thefootnote}{}

\renewcommand{\PaperNumber}{046}

\FirstPageHeading

\ShortArticleName{Companions to Andrews--Gordon--Bressoud Identities and CMPP Conjectures}

\ArticleName{Companions to the Andrews--Gordon\\ and Andrews--Bressoud Identities\\ and Recent Conjectures of Capparelli, Meurman,\\ Primc, and Primc\footnote{This paper is a~contribution to the Special Issue on Recent Advances in Vertex Operator Algebras in honor of James Lepowsky. The~full collection is available at \href{https://sigma-journal.com/Lepowsky.html}{https://sigma-journal.com/Lepowsky.html}}}

\Author{Matthew C. RUSSELL}

\AuthorNameForHeading{M.~C.~Russell}

\Address{Department of Statistics, Texas A\&M University, USA}
\Email{\mail{matthewcrussell@tamu.edu}}
\URLaddress{\url{http://matthewcrussell.github.io}}

\ArticleDates{Received September 05, 2025, in final form April 28, 2026; Published online May 12, 2026}

\Abstract{We find bivariate generating functions for the $k=1$ cases of recently conjectured colored partition identities of Capparelli, Meurman, A.~Primc, and M.~Primc that are slight variants of the generating functions for the sum sides of the Andrews--Gordon and Andrews--Bressoud identities, relating them to recent work of Warnaar. This $k=1$ cases turn out to be equivalent to identities of Jing, Misra, and Savage. Finally, we provide bijections for these identities involving two-rowed cylindric partitions, in the spirit of Corteel.}


\Classification{11P84; 05A17; 17B67}

\begin{flushright}
\begin{minipage}{55mm}
\it To James Lepowsky, on the \\ occasion of his eightieth birthday
\end{minipage}
\end{flushright}

\renewcommand{\thefootnote}{\arabic{footnote}}
\setcounter{footnote}{0}

\section{Introduction}

Let $n$ be a non-negative integer. A {\it partition} of $n$ is a list of integers $\bigl(\lambda_1,\lambda_2,\dots,\lambda_k\bigr)$ such that $ \lambda_1+\lambda_2+\cdots+\lambda_k = n$ and $\lambda_1 \ge \lambda_2 \ge \cdots \ge \lambda_k \ge 1$. We will frequently write this as $\lambda_1+\lambda_2+\cdots+\lambda_k$. For a given partition $\lambda$, let $f_j$ be the number of parts of $\lambda$ equal to $j$. We~will use the term ``weight'' to mean the sum of the parts of a partition or cylindric partition (to be defined later).

Our starting point is the following well-known partition theorem of Basil Gordon~\cite{Gord}, and its analytic form due to George E.~Andrews~\cite{And66} (see also Chapter~7 of Andrews's text~\cite{And-book} or Chapter~3 of Andrew V.~Sills's text~\cite{Sil-book}).

\begin{Theorem} \label{thm:Gordon}
Fix $\ell\ge 1$, and let $0\le i\le \ell$. Let $B_{\ell,i}(n)$ denote the number of partitions of $n$ satisfying, for all $j$, $f_j + f_{j+1}\le \ell$ and $f_1 \le i$. Let $A_{\ell,i}(n)$ denote the number of partitions of $n$ into parts $\not\equiv 0, \pm (i+1)$ $(\bmod~2\ell+3)$. Then $A_{\ell,i}(n)=B_{\ell,i}(n)$ for all $n$.
\end{Theorem}

\begin{Theorem} \label{thm:AndGor}
For $0\le i \le \ell,$
\begin{equation*}
\sum_{n_1,n_2,\dots,n_{\ell}\ge 0} \frac{q^{N_1^2+N_2^2+\cdots+N_{\ell}^2+N_{i+1}+N_{i+2}+\cdots+N_{\ell}}}{(q;q)_{n_1}(q;q)_{n_2}\cdots(q;q)_{n_{\ell}}} = \frac{\bigl(q^{i+1},q^{2\ell+2-i},q^{2\ell+3};q^{2\ell+3}\bigr)_\infty}{(q;q)_\infty}.
\end{equation*}
\end{Theorem}
Here (and throughout the paper), $N_i = n_i +n_{i+1}+\cdots+n_{\ell}$. For the purposes of our present paper, we will view such functions as formal power series in $q$. We are using the standard $q$-Pochhammer notation
$(a;q)_n = \prod_{j=1}^n \bigl(1-aq^{j-1}\bigr)$ for $n \in \NN \cup \{\infty\}$ and $(a_1,\dots,a_m;q)_n=(a_1;q)_n\cdots(a_m;q)_n$. We also define $(a;q)_{-1}=1$.

The generating function in Theorem~\ref{thm:AndGor} was further refined by Andrews~\cite{AndPNAS} into a bivariate generating function. Let $B_{\ell,i,m}(n)$ denote the number of partitions of $n$ with exactly $m$ parts satisfying, for all $j$, $f_j + f_{j+1}\le \ell$ and $f_1 \le i$. Then
\begin{equation} \label{eq:AndGorxq}
\sum_{m,n \ge 0} B_{\ell,i,m}(n) x^m q^n =\sum_{n_1,n_2,\dots,n_{\ell}\ge 0} \frac{x^{N_1+N_2+\cdots+N_{\ell}}q^{N_1^2+N_2^2+\cdots+N_{\ell}^2+N_{i+1}+N_{i+2}+\cdots+N_{\ell}}}{(q;q)_{n_1}(q;q)_{n_2}\cdots(q;q)_{n_{\ell}}}.
\end{equation}
Note that these identities all have odd modulus. The even modulus companions were found by David M. Bressoud~\cite{Bress0,Bress1,Bress2}; these are known as the Andrews--Bressoud identities.
\begin{Theorem}
Fix $\ell\ge 1$, and let $0\le i\le \ell$. Let $B^\star_{\ell,i}(n)$ denote the number of partitions of $n$ satisfying, for all $j$, $f_j + f_{j+1}\le \ell$ and, if $f_j + f_{j+1}=\ell$, then $j\cdot f_j+(j+1)\cdot f_{j+1} \equiv i$ $(\bmod~2)$, and furthermore $f_1 \le i$. Let $A^\star_{\ell,i}(n)$ be the coefficient of $q^n$ in the product \[\frac{\bigl(q^{i+1},q^{2\ell-i+1},q^{2\ell+2};q^{2\ell+2}\bigr)_\infty}{(q;q)_\infty}.\] $($For $i\ne \ell$, this is the number of partitions of $n$ into parts $\not\equiv 0, \pm (i+1)$ $(\bmod~2\ell+2)$.$)$ Then $A^\star_{\ell,i}(n)=B^\star_{\ell,i}(n)$ for all~$n$.
\end{Theorem}

\begin{Theorem} \label{thm:AndBre}
For $0\le i \le \ell,$
\begin{equation*}
\sum_{n_1,n_2,\dots,n_{\ell}\ge 0} \frac{q^{N_1^2+N_2^2+\cdots+N_{\ell}^2+N_{i+1}+N_{i+2}+\cdots+N_{\ell}}}{(q;q)_{n_1}(q;q)_{n_2}\cdots(q;q)_{n_{\ell-1}}\bigl(q^2;q^2\bigr)_{n_{\ell}}} = \frac{\bigl(q^{i+1},q^{2\ell+1-i},q^{2\ell+2};q^{2\ell+2}\bigr)_\infty}{(q;q)_\infty}.
\end{equation*}
\end{Theorem}

Again, this generating function in Theorem~\ref{thm:AndBre} can be refined into a bivariate generating function. Let $B^\star_{\ell,i,m}(n)$ denote the number of partitions of $n$ counted by $B^\star_{\ell,i}(n)$ with exactly $m$ parts. Then
\begin{equation} \label{eq:AndBrexq}
\sum_{m,n \ge 0} B^\star_{\ell,i,m}(n) x^m q^n =\sum_{n_1,n_2,\dots,n_{\ell}\ge 0} \frac{x^{N_1+N_2+\cdots+N_{\ell}}q^{N_1^2+N_2^2+\cdots+N_{\ell}^2+N_{i+1}+N_{i+2}+\cdots+N_{\ell}}}{(q;q)_{n_1}(q;q)_{n_2}\cdots(q;q)_{n_{\ell-1}}\bigl(q^2;q^2\bigr)_{n_{\ell}}}.
\end{equation}
For some recent papers on other facets of the Andrews--Gordon and Andrews--Bressoud identities, see~\cite{ADJM,KLRS,KurParity,Stant,Wiecz}.

We now turn our attention to a recent article~\cite{CMPP} of Stefano Capparelli, Arne Meurman, Andrej Primc, and Mirko Primc, in which they conjectured intriguing (colored) partition identi\-ties related to standard representations of the affine Lie algebra of type \smash{$C_\ell^{(1)}$} (for $\ell\ge 2$). These (conjectured) identities follow in a long vein of research connecting partition identities to the representation theory of affine Lie algebras~\cite{Capp,KanRus-idf,KanRus-stair,LepWil-struI,MP87,MP,Nan-thesis,TakTsu-nandi}. Most directly, these conjectures built off of the work of M.~Primc and Tomislav \v{S}iki\'{c}~\cite{PS2,PS1} and Goran Trup\v{c}evi\'{c}~\cite{Tru}. M.~Primc and Trup\v{c}evi\'{c} have very recently proven the type \smash{$C_\ell^{(1)}$} conjectures in the special case of initial conditions $[k,0,\dots,0]$~\cite{PT}.

First, let us consider the conjectures in Section~4 of their paper~\cite{CMPP}. Capparelli, Meurman, A.~Primc, and M.~Primc left the representation-theoretic interpretations of this family of conjectures as an open question; Shashank Kanade, S. Ole Warnaar, Shunsuke Tsuchioka, and the author~\cite{KRTW} speculate that these conjectures have an association with type \smash{$A_{2n}^{(2)}$}. Fix $\ell\ge 2$, and construct the following array with $2\ell$ rows (and $\ell$ copies of $\NN$):
\begin{equation*}
\begin{matrix}
1 && 3 && 5 && 7 && 9 && 11 & \\
& 2 && 4 && 6 && 8 && 10 && 12\\
1 && 3 && 5 && 7 && 9 && 11 & \\
& 2 && 4 && 6 && 8 && 10 && 12\\
1 && 3 && 5 && 7 && 9 && 11 & \\
& 2 && 4 && 6 && 8 && 10 && 12\\
\vdots & \vdots &\vdots& \vdots &\vdots& \vdots &\vdots& \vdots &\vdots& \vdots &\vdots& \vdots\\
1 && 3 && 5 && 7 && 9 && 11 & \\
& 2 && 4 && 6 && 8 && 10 && 12\\
1 && 3 && 5 && 7 && 9 && 11 & \\
& 2 && 4 && 6 && 8 && 10 && 12
 \end{matrix}\quad\dots\,.
\end{equation*}
Parts with the same value but lying in different rows are distinct, and can be thought of as having distinct colors. We will sometimes call these colored partitions.
To the above array, we will associate an extended array of frequencies
\begin{equation*}
\begin{matrix}
k_1 && f_1 && f_3 && f_5 && f_7 && f_9 && f_{11} & \\
& \cdot && f_2 && f_4 && f_6 && f_8 && f_{10} && f_{12}\\
k_3 &&f_1 && f_3 && f_5 && f_7 && f_9 && f_{11} & \\
& \cdot && f_2 && f_4 && f_6 && f_8 && f_{10} && f_{12}\\
k_5 &&f_1 && f_3 && f_5 && f_7 && f_9 && f_{11} & \\
& \cdot && f_2 && f_4 && f_6 && f_8 && f_{10} && f_{12}\\
\vdots & \vdots &\vdots& \vdots &\vdots& \vdots &\vdots& \vdots &\vdots& \vdots &\vdots& \vdots&\vdots& \vdots\\
k_4 && f_1 && f_3 && f_5 && f_7 && f_9 && f_{11} & \\
& \cdot && f_2 && f_4 && f_6 && f_8 && f_{10} && f_{12} \\
k_2 && f_1 && f_3 && f_5 && f_7 && f_9 && f_{11} & \\
& k_0 && f_2 && f_4 && f_6 && f_8 && f_{10} && f_{12}
 \end{matrix}\quad\dots,
\end{equation*}
where each $f_j$ indicates how many times the corresponding part occurs in the original array. The nonnegative integers $k_0,\dots,k_\ell$ provide initial conditions. For odd $i\le \ell$, $k_i$ is in the $i$-th row, while for even $0<i \le \ell$, $k_i$ is in the $(2\ell-i+1)$-th row, and $k_0$ is always in the $(2\ell)$-th row. Note that our ordering of $k_0,\dots,k_\ell$ is different than that of the original paper. The symbol~$\cdot$~is a~placeholder that represents parts that are forbidden from appearing. Every entry in the second column, other than the final one, has a $\cdot$ in it. Here is the array in the case $\ell=4$:
\begin{equation*}
\begin{matrix}
k_1 && f_1 && f_3 && f_5 && f_7 && f_9 && f_{11} & \\
& \cdot && f_2 && f_4 && f_6 && f_8 && f_{10} && f_{12}\\
k_3 &&f_1 && f_3 && f_5 && f_7 && f_9 && f_{11} & \\
& \cdot && f_2 && f_4 && f_6 && f_8 && f_{10} && f_{12}\\
k_4 && f_1 && f_3 && f_5 && f_7 && f_9 && f_{11} & \\
& \cdot && f_2 && f_4 && f_6 && f_8 && f_{10} && f_{12} \\
k_2 && f_1 && f_3 && f_5 && f_7 && f_9 && f_{11} & \\
& k_0 && f_2 && f_4 && f_6 && f_8 && f_{10} && f_{12}
 \end{matrix}\quad\dots\,.
\end{equation*}
A downward path $\mathcal Z$ in the extended frequency array is a set of frequencies containing exactly one element in each row, where the elements in consecutive rows are adjacent to each other. (The placeholders $\cdot$ and the initial conditions $k_i$ are allowed to be part of these downward paths.) Here are two sample downward paths in the case with $\ell=3$:
\begin{equation*}
\begin{matrix}
\color{blue}{\bf k_1} && f_1 && {f_3} && f_5 && \color{red}{\bf f_7} && f_9 && f_{11} & \\
& \color{blue}{\bf \bullet} && f_2 && f_4 && f_6 && \color{red}{\bf f_8} && f_{10} && f_{12}\\
k_3 &&\color{blue}{\bf f_1} && {f_3} && f_5 && f_7 && \color{red}{\bf f_9} && f_{11} & \\
& \cdot && \color{blue}{\bf f_2} && f_4 && f_6 && f_8 && \color{red}{\bf f_{10}} && f_{12}\\
k_2 && \color{blue}{\bf f_1} && f_3 && f_5 && f_7 && \color{red}{\bf f_9} && f_{11} & \\
& k_0 && \color{blue}{\bf f_2} && f_4 && f_6 && \color{red}{\bf f_8} && f_{10} && f_{12}
 \end{matrix}\quad\dots\,.
\end{equation*}
 We say that an array of frequencies is $[k_0,\dots,k_\ell]$-admissible if, for all downward paths $\mathcal Z$ in the extended frequency array,
\begin{equation*}
\sum_{m \in \mathcal Z} m \le k,
\end{equation*}
where $k = \sum_{i=0}^\ell k_i$.

Fix $k=1$. For $0 \le i \le \ell$, let $F(i,j,n)$ be the number of $\left[k_0,k_1,\dots,k_\ell\right]=\left[0,\dots,0,1,0,\dots,0\right]$-admissible colored partitions of $n$ with exactly $j$ parts, where $k_i=1$ (and all others are 0). Then, define $P_{i}(z,q)$ to be the bivariate generating function
\begin{equation*}
P_{i}(z,q)=\sum_{n,j\ge 0} F(i,j,n)z^j q^n.
\end{equation*}
The first goal of this paper will be to prove the following multisum for $P_i(z,q)$.
\begin{Theorem} \label{thm:mainresult1}
For $0\le i \le \ell,$
\begin{equation} \label{eq:mainresult1}
P_i(z,q)=\sum_{n_1,n_2,\dots,n_{\ell}\ge 0} \frac{z^{N_1}q^{N_1^2+N_2^2+\cdots+N_{\ell}^2+N_{i+1}+N_{i+2}+\cdots+N_{\ell}}}{(q;q)_{n_1}(q;q)_{n_2}\cdots(q;q)_{n_{\ell}}}.
\end{equation}
\end{Theorem}

Take a moment to appreciate just how similar the bivariate multisum in Theorem~\ref{thm:mainresult1} is to the multisum of~\eqref{eq:AndGorxq}. The only difference is that the factor $x^{N_1+N_2+\cdots+N_{\ell}}$ has changed to $z^{N_1}$. As an immediate corollary,
\begin{align}
P_i(1,q)& =\sum_{n_1,n_2,\dots,n_{\ell}\ge 0} \frac{q^{N_1^2+N_2^2+\cdots+N_{\ell}^2+N_{i+1}+N_{i+2}+\cdots+N_{\ell}}}{(q;q)_{n_1}(q;q)_{n_2}\cdots(q;q)_{n_{\ell}}}\nonumber\\
& = \frac{\bigl(q^{i+1},q^{2\ell+2-i},q^{2\ell+3};q^{2\ell+3}\bigr)_\infty}{(q;q)_\infty}\label{eq:sumprod1}
\end{align}
by Theorem~\ref{thm:AndGor}.

Now, let us consider the conjectured identities in Section~3 of Capparelli, Meurman, A.~Primc, and M.~Primc~\cite{CMPP}, along with its sequel by M.~Primc~\cite{primc2023}. We still fix $\ell\ge 2$, but now construct arrays with $2\ell-1$ rows. Section~3 of~\cite{CMPP} dealt with arrays whose top and bottom rows consist of odd integers:
\begin{equation*}
\begin{matrix}
1 && 3 && 5 && 7 && 9 && 11 & \\
& 2 && 4 && 6 && 8 && 10 && 12\\
1 && 3 && 5 && 7 && 9 && 11 & \\
& 2 && 4 && 6 && 8 && 10 && 12\\
\vdots & \vdots &\vdots& \vdots &\vdots& \vdots &\vdots& \vdots &\vdots& \vdots &\vdots& \vdots\\
1 && 3 && 5 && 7 && 9 && 11 & \\
& 2 && 4 && 6 && 8 && 10 && 12 \\
1 && 3 && 5 && 7 && 9 && 11 &
 \end{matrix}\quad\dots\,.
\end{equation*}

In~\cite{primc2023}, the top and bottom rows consist of even integers:
\begin{equation*}
\begin{matrix}
& 2 && 4 && 6 && 8 && 10 && 12\\
1 && 3 && 5 && 7 && 9 && 11 & \\
& 2 && 4 && 6 && 8 && 10 && 12\\
\vdots & \vdots &\vdots& \vdots &\vdots& \vdots &\vdots& \vdots &\vdots& \vdots &\vdots& \vdots\\
& 2 && 4 && 6 && 8 && 10 && 12 \\
1 && 3 && 5 && 7 && 9 && 11 & \\
& 2 && 4 && 6 && 8 && 10 && 12
 \end{matrix}\quad\dots\,.
\end{equation*}

The $k=1$ case of the conjectures of Capparelli, Meurman, A.~Primc, and M.~Primc (corresponding to odd first and last rows) was proved by Jehanne Dousse and Isaac Konan~\cite{DouKon} using perfect crystals. This is the same case that we will prove, but our proof methods are different, and we will simultaneously prove the two cases (odd first and last rows and even first and last rows).

We again form arrays of frequencies, and attach nonnegative integers $k_0,\dots,k_\ell$. Our definition of an admissible partition, using downward paths, is the same as above, and we once again restrict our attention to the case $k=1$. See Section~\ref{sec:partitions4} for specific details, including directions for how $k_0,\dots,k_\ell$ are attached. Once again, for $0 \le i \le \ell$, let $F^\star(i,j,n)$ be the number of $\left[k_0,k_1\dots,k_\ell\right]=\left[0,\dots,0,1,0,\dots,0\right]$-admissible colored partitions of $n$ with exactly $j$ parts, where the 1 in $\left[0,\dots,0,1,0,\dots,0\right]$ is located in position $i$. Then, define $P^\star_{i}(z,q)$ to be the bivariate generating function
\begin{equation*}
P^\star_{i}(z,q)=\sum_{n,j\ge 0} F^\star(i,j,n)z^j q^n.
\end{equation*}
We then have the corresponding result here for $P^\star_i(z,q)$.
\begin{Theorem} \label{thm:mainresult2}
For $0\le i \le \ell$,
\begin{equation} \label{eq:mainresult2}
P^\star_i(z,q)=\sum_{n_1,n_2,\dots,n_{\ell}\ge 0} \frac{z^{N_1}q^{N_1^2+N_2^2+\cdots+N_{\ell}^2+N_{i+1}+N_{i+2}+\cdots+N_{\ell}}}{(q;q)_{n_1}(q;q)_{n_2}\cdots(q;q)_{n_{\ell-1}}\bigl(q^2;q^2\bigr)_{n_{\ell}}}.
\end{equation}
\end{Theorem}
This again is very similar to~\eqref{eq:AndBrexq} and gives us the sum-to-product identity
\begin{align}
P^\star_i(1,q)& =\sum_{n_1,n_2,\dots,n_{\ell}\ge 0} \frac{q^{N_1^2+N_2^2+\cdots+N_{\ell}^2+N_{i+1}+N_{i+2}+\cdots+N_{\ell}}}{(q;q)_{n_1}(q;q)_{n_2}
\cdots(q;q)_{n_{\ell-1}}\bigl(q^2;q^2\bigr)_{n_{\ell}}}\nonumber\\
& = \frac{\bigl(q^{i+1},q^{2\ell+1-i},q^{2\ell+2};q^{2\ell+2}\bigr)_\infty}{(q;q)_\infty}\label{eq:sumprod2}
\end{align}
by Theorem~\ref{thm:AndBre}.

Actually, the $k=1$ cases of both families of conjectures turn out to be equivalent to a~theorem of Naihuan Jing, Kailash C. Misra, and Carla D. Savage published in 2001~\cite{JMS}. This fact appears to have gone unnoticed for quite some time, as the notation is a bit different in the papers (Jing, Misra, and Savage's difference conditions are expressed using absolute values of differences of colors, not in terms of the ``downward paths'' of Capparelli, Meurman, A.~Primc, and M.~Primc). The author is indebted to Shunsuke Tsuchioka for pointing this out to him after the first preprint of this paper appeared. Thus, this present paper cannot take credit for initially proving~\eqref{eq:sumprod1} and~\eqref{eq:sumprod2}, as those results first appeared in~\cite{JMS}. However, Theorems~\ref{thm:mainresult1} and~\ref{thm:mainresult2} (the bivariate multisums) are new. We summarize which cases of these Capparelli--Meurman--Primc--Primc conjectures have been proven in Table~\ref{table} (see Section~\ref{sec:partitions4} for more information about the ``odd-odd'' frequency arrays).

\begin{table}[ht] \label{table}
 \centering\small\renewcommand{\arraystretch}{1.2}
 \begin{tabular}{|c|c|} \hline
all $k$, 3 rows, odd-odd arrays & Meurman and Primc~\cite{MP} \\ \hline
$k=1$, all cases & Jing, Misra, and Savage~\cite{JMS} \\ \hline
$\bigl(k_1,k_2,\dots,k_\ell\bigr)=(1, 0,\dots , 0)$, $2\ell-1$ rows&Primc and \v{S}iki\'{c}~\cite{PS2} \\ \hline
$k=1$, $2\ell-1$ rows, odd-odd arrays & Dousse and Konan~\cite{DouKon} \\ \hline
$\bigl(k_1,k_2,\dots,k_\ell\bigr)=(k, 0,\dots , 0)$, cases from~\cite{CMPP}&Primc and Trup\v{c}evi\'{c}~\cite{PT} \\\hline
$k=1$, all cases & This paper \\\hline
 \end{tabular}
 \caption{A summary of proven cases of the Capparelli--Meurman--Primc--Primc conjectures.}
\end{table}

The final ingredients for this paper are cylindric partitions. Cylindric partitions were first introduced by Ira Gessel and Christian Krattenthaler~\cite{GesKra}. For a given composition $(c_1,c_2,\dots,c_r)$, where each $c_j$ is a nonnegative integer, an $r$-rowed cylindric partition of profile $(c_1,c_2,\dots,c_r)$ is a list of $r$ partitions \smash{$\bigl(\tau^{(1)},\tau^{(2)},\dots,\tau^{(r)}\bigr)$} satisfying \smash{$\tau^{(j)}_i\geq \tau^{(j+1)}_{i+c_{(j+1)}}$} for $1\le j <r$ and all $i$, along with \smash{$\tau^{(r)}_i\geq \tau^{(1)}_{i+c_1}$} for all $i$. In this paper, we concern ourselves only with two-rowed cylindric partitions. In this case, the definition simplifies to a two-rowed cylindric partition $\Lambda$ of profile $(c_1,c_2)$, where $c_1,c_2\ge0$, is a pair of partitions $\bigl(\tau^{(1)},\tau^{(2)}\bigr)$ such that, for all $i$,
\begin{equation*}
\tau_i^{(1)} \ge \tau_{i+c_2}^{(2)} \qquad \text{and} \qquad \tau_i^{(2)} \ge \tau_{i+c_1}^{(1)}.
\end{equation*}
For example, if $\tau^{(1)}=6+3+2+2$ and $\tau^{(2)}=5+2+1$, then $\bigl(\tau^{(1)},\tau^{(2)}\bigr)$ is a cylindric partition of profile $(2,1)$. We can see this by arranging the parts in the following way:
\begin{equation*}
\begin{matrix}
\textcolor{gray}{\tau^{(2)} \rightarrow } & \textcolor{gray}{\text{ }}& \textcolor{gray}{\text{ }}& \textcolor{gray}{\text{ }}&\textcolor{gray}{5} & \textcolor{gray}{2} & \textcolor{gray}{1} \\
\tau^{(1)} \rightarrow & \text{ }& 6 & 3 & 2 &2 \\
\tau^{(2)} \rightarrow & 5& 2 & 1
\end{matrix}
\end{equation*}
This satisfies the property that the parts are arranged in a nonincreasing fashion along each row and each column.
Note that we repeat $\tau^{(2)}$ in grey above $\tau^{(1)}$ to demonstrate the inequalities arising from \smash{$\tau_i^{(2)} \ge \tau_{i+c_1}^{(1)}$}.

Alexei Borodin~\cite{Bor} found infinite product representations for the (single-variable) generating functions of cylindric partitions with a given profile. In the case of two-rowed cylindric partitions, this ends up as $(q;q)^{-1}_\infty$ times the generating function for either the Andrews--Gordon or Andrews--Bressoud identities, depending on the parity of $c_1+c_2$. Omar Foda and Trevor Welsh~\cite{FodWel} then gave a proof of the Andrews--Gordon identities explicitly using two-rowed cylindric partitions, and noted that the infinite product forms of the generating functions for cylindric partitions of a given profile can be deduced using the work of Gessel and Krattenthaler~\cite{GesKra}, along with the appropriate Macdonald identity. Sylvie Corteel and Trevor Welsh~\cite{CW} deduced functional equations governing bivariate generating functions for cylindric partitions. In recent work~\cite{WarnCyl}, Warnaar found finitizations of bivariate multisum generating functions for two-rowed cylindric partitions. For these bivariate $(z,q)$-generating functions, the exponent of $z$ keeps track of the size of the largest part of the cylindric partition (while the exponent of $q$, of course, counts the sum of all of the parts).

As it turns out, \[ \frac{P_i(z,q)}{(zq;q)_\infty} \qquad \text{and} \qquad \displaystyle \frac{P_i^\star (z,q)}{(zq;q)_\infty}\] turn out to be exactly the same as the bivariate generating functions found by Warnaar~\cite{WarnCyl}: the multisums given in Theorems~\ref{thm:mainresult1} and~\ref{thm:mainresult2} can be obtained by taking the limit as $L\to \infty$ in equations (7.24) and (7.25) of his paper.

Partition bijections are always treasured in the field; we note the bijective proofs of the Rogers--Ramanujan identities offered by Adriano Garsia and Stephen Milne~\cite{GarMil} and Bressoud and Doron Zeilberger~\cite{BZ}. Igor Pak~\cite{PakSurvey} offers a comprehensive survey of partition bijections. Remarkably, Corteel~\cite{Cort} gave a bijective proof of the second Rogers--Ramanujan identity divided by $(q;q)_\infty$:
\begin{equation*}
\frac{1}{(q;q)_\infty}\sum_{n=0}\frac{q^{n^2+n}}{(q;q)_n} = \frac{1}{(q;q)_\infty\bigl(q^2,q^3;q^5\bigr)_\infty}
\end{equation*}
using two-rowed cylindric partitions of profile $(3,0)$. Inspired by her work, the last result of our paper is a bijection between two-rowed cylindric partitions of profile $(2\ell+1,0)$ and pairs $(\lambda,\mu)$, where $\lambda$ is an ordinary partition and $\mu$ is a $[k_0,k_1,\dots,k_\ell]=[1,0,\dots,0]$-admissible partition of Capparelli, Meurman, Primc, and Primc on $2\ell$ rows.

The rest of the paper is organized as follows. In Section~\ref{sec:partitions2}, we produce functional equations for the generating functions for the admissible partitions on arrays with $2\ell$ rows, corresponding to the Andrews--Gordon case. In Section~\ref{sec:comp}, we complete our proof of Theorem \ref{thm:mainresult1} by arguing that the appropriate bivariate multisums do, in fact, satisfy the same functional equations as the generating functions from Section~\ref{sec:partitions2}. Sections~\ref{sec:partitions4} and \ref{sec:comp5} repeat this work for arrays on $2\ell-1$ rows (corresponding to the Andrews--Bressoud case). Section~\ref{sec:6} deals with two-rowed cylindric partitions and the aforementioned bijective arguments. Section~\ref{sec:conc7} provides some thoughts about further research suggested by this paper. Appendix~\ref{sec:appendix} illustrates several of the computations of this paper by explicitly demonstrating them in the case $\ell=4$.

\subsection{Jim Lepowsky's influence on the author's mathematical journey}

It is only fitting that this introduction should conclude with a note of thanks to Jim Lepowsky, to whom this paper is dedicated. The author is extremely grateful to Jim for teaching a course in the spring of 2014 on the theory of partitions and vertex operator algebras. Jim's original plan for the course was to spend roughly the first half on integer partitions and the second half on VOAs. As it turned out, out of the six students in the course, the author was the only one who was unfamiliar with VOAs. Very wisely, Jim changed gears for the course at the halfway point and organized a class project on a motivated proof of the G\"ollnitz--Gordon--Andrews identities that was later published~\cite{CKLMQRS}. During that same semester, the author began working on two more projects that were ultimately published: one joint with Jim, Shashank Kanade, and Drew Sills~\cite{KLRS}, and one joint with Kanade that reflected lots of discussion with Jim~\cite{KanRus-idf}. So, the author's first three publications on integer partitions are all directly attributable to that one course of Jim's. Thank you!

\section{Partitions: Andrews--Gordon case}\label{sec:partitions2}

Fix $\ell$ and take an array with $2\ell$ rows to it (with $\ell$ copies of $\NN$). Define $P_i(z,q)$ to be the bivariate generating function for $[k_0,k_1,\dots,k_\ell]=[0,\dots,0,1,0,\dots,0]$-admissible partitions, where below $k_i=1$ and all other $k_j=0$:
\begin{equation*}
\begin{matrix}
k_1 && f_1 && f_3 && f_5 && f_7 && f_9 && f_{11} & \\
& \cdot && f_2 && f_4 && f_6 && f_8 && f_{10} && f_{12}\\
k_3 &&f_1 && f_3 && f_5 && f_7 && f_9 && f_{11} & \\
& \cdot && f_2 && f_4 && f_6 && f_8 && f_{10} && f_{12}\\
\vdots & \vdots &\vdots& \vdots &\vdots& \vdots &\vdots& \vdots &\vdots& \vdots &\vdots& \vdots&\vdots& \vdots\\
k_2 && f_1 && f_3 && f_5 && f_7 && f_9 && f_{11} & \\
& k_0 && f_2 && f_4 && f_6 && f_8 && f_{10} && f_{12}
 \end{matrix}\quad\dots\,.
\end{equation*}
Here, for $i$ odd, $k_i$ occurs in the $i$-th row, while for $i$ even and $i\ge 2$, $k_i$ occurs in the $(2\ell+1-i)$-th row. Finally, $k_0$ always occurs in the $(2\ell)$-th row. Note that, for $1 \le i \le 2\ell$, the only allowable occurrence of the part 1 is in row $i$ (if $i$ is odd) or row $2\ell+1-i$ (if $i$ is even).

We will now show that these satisfy the following system of functional equations.
\begin{Theorem} \label{thm:funceq1A}
The generating functions $P_i(z,q)$ satisfy the following functional equations:
\begin{align}
&P_i(z,q) - P_{i+1}(zq,q)= \sum_{j=1}^i zq^j P_{i-j+1}\bigl(zq^{j+1},q\bigr), \qquad 0 \le i \le \ell-1, \label{eq:origFE_gf1} \\
&P_\ell(z,q) - P_\ell(zq,q)= \sum_{j=1}^\ell zq^j P_{\ell-j+1}\bigl(zq^{j+1},q\bigr). \label{eq:origFE_gf2}
\end{align}
\end{Theorem}
Throughout the paper, empty sums are defined to equal zero. As an example, the instance of \eqref{eq:origFE_gf1} for $i=0$ simply states that $P_0(z,q)-P_1(zq,q)=0$. An illustration of the following proof of this theorem for $\ell=4$ is given in Appendix~\ref{subsec:combspec}.

\begin{proof}
First, consider the case $i=0$. We will show that $P_0(z,q)=P_1(zq,q)$ by means of a~``flipping'' trick. We present the extended frequency arrays for $P_0(z,q)$ and $P_1(z,q)$, with their initial conditions shown as $\ficone$, and $\cdot$ representing parts that are forbidden from appearing:
\begin{align*}
&P_0(z,q)\colon \quad
\begin{matrix}
\vdots & \vdots &\vdots& \vdots &\vdots& \vdots &\vdots& \vdots &\vdots& \vdots &\vdots& \vdots&\vdots& \vdots\\
& \cdot && \cdot && \cdot && f_6 && f_8 && f_{10} && f_{12} \\
\cdot && \cdot && \cdot && f_5 && f_7 && f_9 && f_{11} & \\
& \cdot && \cdot && f_4 && f_6 && f_8 && f_{10} && f_{12} \\
\cdot && \cdot && f_3 && f_5 && f_7 && f_9 && f_{11} & \\
& \ficone && f_2 && f_4 && f_6 && f_8 && f_{10} && f_{12}\vspace{3mm}
 \end{matrix}\quad\dots\,, \\
&P_1(z,q)\colon \quad
\begin{matrix}
\ficone && f_1 && f_3 && f_5 && f_7 && f_9 && f_{11} & \\
& \cdot && f_2 && f_4 && f_6 && f_8 && f_{10} && f_{12}\\
\cdot &&\cdot && f_3 && f_5 && f_7 && f_9 && f_{11} & \\
& \cdot && \cdot && f_4 && f_6 && f_8 && f_{10} && f_{12}\\
\cdot && \cdot && \cdot && f_5 && f_7 && f_9 && f_{11} & \\
& \cdot && \cdot && \cdot && f_6 && f_8 && f_{10} && f_{12}\\
\vdots & \vdots &\vdots& \vdots &\vdots& \vdots &\vdots& \vdots &\vdots& \vdots &\vdots& \vdots&\vdots& \vdots\\ \end{matrix}\quad\dots\,.
\end{align*}
Now, $P_1(zq,q)$ is the generating function for partitions counted by $P_1(z,q)$ with 1 added to every part:
\begin{equation*}
P_1(zq,q)\colon \quad
\begin{matrix}
\ficone && f_2 && f_4 && f_6 && f_8 && f_{10} && f_{12} & \\
& \cdot && f_3 && f_5 && f_7 && f_9 && f_{11} && f_{13}\\
\cdot &&\cdot && f_4 && f_6 && f_8 && f_{10} && f_{12} & \\
& \cdot && \cdot && f_5 && f_7 && f_9 && f_{11} && f_{13}\\
\cdot && \cdot && \cdot && f_6 && f_8 && f_{10} && f_{12} & \\
& \cdot && \cdot && \cdot && f_7 && f_9 && f_{11} && f_{13}\\
\vdots & \vdots &\vdots& \vdots &\vdots& \vdots &\vdots& \vdots &\vdots& \vdots &\vdots& \vdots&\vdots& \vdots\\ \end{matrix}\quad\dots\,.
\end{equation*}
And so flipping this picture upside down shows that $P_0(z,q)$ and $P_1(zq,q)$ count the same partitions, and thus $P_0(z,q)=P_1(zq,q)$. As noted right before the beginning of the proof, the $i=0$ instance of~\eqref{eq:origFE_gf1} simplifies to $P_0(z,q)-P_1(zq,q)=0$, so we are done with the $i=0$ case.

Next, consider the case where $i$ is odd (and $i\ne \ell$). Suppose that we have an admissible partition counted by $P_i(z,q)$. Consider, for $1 \le j \le i$, the set of parts $j$ that may occur in row $i-j+1$; at most one of these parts can occur in this partition (as they are all in the same downward path). Here is a picture corresponding to the partitions counted by $P_i(z,q)$; the parts mentioned in the previous sentence are shown in blue. Recall that $k_i$ (represented by the $\ficone$) occurs in row $i$:
\begin{gather*}
P_i(z,q)\colon \quad
\begin{matrix}
\!\cdots\! && \!\cdots\! && \!\cdots\! && \cdot && \color{blue}{f_{i}} && f_{i+2} && f_{i+4} & \\
& \!\cdots\! && \!\cdots\! && \cdot && \color{blue}{f_{i-1}} && f_{i+1} && f_{i+3} && f_{i+5}\\
\vdots & \vdots &\vdots& \vdots &\vdots& \vdots &\iddots& \vdots &\vdots& \vdots &\vdots& \vdots&\vdots& \vdots\\
& \cdot && \cdot && \color{blue}{f_4} && f_6 && f_8 && f_{10} && f_{12}\\
\cdot &&\cdot && \color{blue}{f_3} && f_5 && f_7 && f_9 && f_{11} & \\
& \cdot && \color{blue}{f_2} && f_4 && f_6 && f_8 && f_{10} && f_{12}\\
\ficone &&\color{blue}{f_1} && f_3 && f_5 && f_7 && f_9 && f_{11} & \\
& \cdot && f_2 && f_4 && f_6 && f_8 && f_{10} && f_{12}\\
\cdot &&\cdot && f_3 && f_5 && f_7 && f_9 && f_{11} & \\
& \cdot && \cdot && f_4 && f_6 && f_8 && f_{10} && f_{12}\\
\vdots & \vdots &\vdots& \vdots &\vdots& \vdots &\vdots& \vdots &\vdots& \vdots &\vdots& \vdots&\vdots& \vdots
 \end{matrix}\quad \dots\,.
\end{gather*}
Suppose that $j$ does occur in the $(i-j+1)$-th row. Here is the corresponding picture, where the red $\color{red}{j}$ indicates the part for which $f_j=1$, while the red $\rs$s show the parts that cannot appear due to the presence of the $\color{red}{j}$ (because they share a downward path with the $\color{red}{j}$):
\begin{gather*}
\begin{matrix}
\cdots && \cdots && \cdots && \cdots && \cdot && \rs && f_{i+2} & \\
& \cdots && \cdots && \cdots && \cdot && \rs && f_{i+1} && f_{i+3} \\
\vdots & \vdots &\vdots& \vdots &\vdots& \vdots &\vdots& \vdots &\iddots& \vdots &\vdots& \vdots&\vdots& \vdots\\
& \cdot && \cdot && \cdot && \rs && f_{j+3} && f_{j+5} && f_{j+7}\\
\cdot &&\cdot && \cdot && \color{red}{j} && f_{j+2} && f_{j+4} && f_{j+6} & \\
& \cdot && \cdot && \rs && \rs && f_{j+3} && f_{j+5} && f_{j+7}\\
\vdots & \vdots &\vdots& \vdots &\iddots& \vdots &\vdots& \vdots &\ddots& \vdots &\vdots& \vdots&\vdots& \vdots\\
& \cdot && \rs && \rs && \rs && \cdots && \cdots && \cdots\\
\ficone &&\rs && \rs && \rs && \rs && \cdots && \cdots & \\
& \cdot && \rs && \rs && \rs && \rs && \cdots && \cdots\\
\cdot &&\cdot && \rs && \rs && \rs && \rs && \cdots & \\
& \cdot && \cdot && \rs && \rs && \rs && \rs &&\cdots\\
\vdots & \vdots &\vdots& \vdots &\vdots& \vdots &\vdots& \vdots &\vdots& \vdots &\vdots& \vdots&\vdots& \vdots\\
 \end{matrix}\quad\dots.
\end{gather*}
After removing the $\color{red}{j}$, the remaining parts will be counted by $P_{i-j+1}\bigl(zq^{j+1},q\bigr)$. To see this, note that the smallest allowable part is in row $i-j+1$. If~$j$ is odd, then $i-j+1$ will also be odd, and we use the previously-observed fact that the only allowable occurrence of the part 1 in the array for $P_{i-j+1}(z,q)$ is in row $i-j+1$. Shifting $z \mapsto zq^{j+1}$ then makes the smallest allowable part equal to~$j+2$. On the other hand, if $j$ is even, then $i-j+1$ will also be even, and we will need to flip the array upside down. After the flip, the smallest allowable part is now in row $2\ell-i+j$. As $i-j+1$ is even, the only allowable occurrence of 1 in the array for $P_{i-j+1}(z,q)$ is in row $2\ell+1-(i-j+1)=2\ell-i+j$. Again, shifting $z \mapsto zq^{j+1}$ then makes the smallest allowable part equal to $j+2$. Thus, if $j$ does occur in the $(i-j+1)$-th row, then, in both cases, that adds a total of $zq^j P_{i-j+1}\bigl(zq^{j+1},q\bigr)$ to the generating function, where the factor of~$zq^j$ comes from the part $j$.

Now, if none of the parts $j$ in the $(i-j+1)$-th rows actually occur, then we have a partition that could be counted by $P_{i+1}(zq,q)$:
\begin{gather*}
P_{i+1}(zq,q)\colon \quad
\begin{matrix}
\!\cdots\! && \!\cdots\! && \!\cdots\! && \!\cdots\! && \cdot && f_{i+2} && f_{i+4} & \\
& \!\cdots\! && \!\cdots\! && \!\cdots\! && \cdot && f_{i+1} && f_{i+3} && f_{i+5}\\
\vdots & \vdots &\vdots& \vdots &\vdots& \vdots &\vdots& \vdots &\vdots& \vdots &\vdots& \vdots&\vdots& \vdots\\
& \cdot && \cdot && \cdot && f_6 && f_8 && f_{10} && f_{12}\\
\cdot &&\cdot && \cdot && f_5 && f_7 && f_9 && f_{11} & \\
& \cdot && \cdot && f_4 && f_6 && f_8 && f_{10} && f_{12}\\
\ficone &&\cdot && f_3 && f_5 && f_7 && f_9 && f_{11} & \\
& \cdot && f_2 && f_4 && f_6 && f_8 && f_{10} && f_{12}\\
\cdot &&\cdot && f_3 && f_5 && f_7 && f_9 && f_{11} & \\
& \cdot && \cdot && f_4 && f_6 && f_8 && f_{10} && f_{12}\\
\vdots & \vdots &\vdots& \vdots &\vdots& \vdots &\vdots& \vdots &\vdots& \vdots &\vdots& \vdots&\vdots& \vdots\\
 \end{matrix}\quad \dots\,.
\end{gather*}
Again, to see this, we use the flipping trick: the smallest allowable part (the part 2) is in row $i+1$, which is an even-indexed row (since $i$ is odd). Flipping moves this part to row $2\ell+1-i$. The partitions counted by $P_{i+1}(z,q)$ only are allowed to have the part 1 in row $2\ell+1-i$; shifting $z\mapsto zq$ increases all parts by 1.

Since exactly one of the above cases ($j$ occurring in row $i-j+1$ for $1 \le j \le i$, or none of those $j$s occurs), we obtain the functional equation $P_i(z,q) = \sum_{j=1}^i zq^j P_{i-j+1}\bigl(zq^{j+1},q\bigr)+ P_{i+1}(zq,q)$.

Next, consider the case where $i$ is even (and $i\not\in\{0,\ell\}$). Suppose that we have an admissible partition counted by $P_i(z,q)$. Consider, for $1 \le j \le i$, the set of parts $j$ that may occur in row $2\ell-i+j$; at most one of these parts can occur in this partition. Here is a picture corresponding to the partitions counted by $P_i(z,q)$; the parts mentioned in the previous sentence are shown in blue. Recall that $k_i$ (represented by the $\ficone$) occurs in row $2\ell+1-i$.
\begin{gather*}
P_i(z,q)\colon \quad
\begin{matrix}
\vdots & \vdots &\vdots& \vdots &\vdots& \vdots &\vdots& \vdots &\vdots& \vdots &\vdots& \vdots&\vdots& \vdots\\
& \cdot && \cdot && f_4 && f_6 && f_8 && f_{10} && f_{12}\\
\cdot &&\cdot && f_3 && f_5 && f_7 && f_9 && f_{11} & \\
& \cdot && {f_2} && f_4 && f_6 && f_8 && f_{10} && f_{12}\\
\ficone &&\color{blue}{f_1} && f_3 && f_5 && f_7 && f_9 && f_{11} & \\
& \cdot && \color{blue}{f_2} && f_4 && f_6 && f_8 && f_{10} && f_{12}\\
\cdot &&\cdot && \color{blue}{f_3} && f_5 && f_7 && f_9 && f_{11} & \\
& \cdot && \cdot && \color{blue}{f_4} && f_6 && f_8 && f_{10} && f_{12}\\
\vdots & \vdots &\vdots& \vdots &\vdots& \vdots &\ddots& \vdots &\vdots& \vdots &\vdots& \vdots&\vdots& \vdots \\
& \!\cdots\! && \!\cdots\! && \cdot && \color{blue}{f_{i-1}} && f_{i+1} && f_{i+3} && f_{i+5}\\
\!\cdots\! && \!\cdots\! && \!\cdots\! && \cdot && \color{blue}{f_{i}} && f_{i+2} && f_{i+4} & \\
 \end{matrix}\quad\dots\,.
\end{gather*}
Suppose that $j$ does occur in the $(2\ell-i+j)$-th row. After removing this $j$, the remaining parts will be counted by $P_{i-j+1}\bigl(zq^{j+1},q\bigr)$. To see this, note that the smallest allowable part is in row $2\ell-i+j$. If $j$ is odd, then $i-j+1$ will be even, and we use the fact that the only allowable occurrence of the part 1 in the array for $P_{i-j+1}(z,q)$ is in row $2\ell+1-(i-j+1)=2\ell-i+j$. Shifting $z \mapsto zq^{j+1}$ then makes the smallest allowable part for this array equal to $j+2$. On the other hand, if $j$ is even, then $i-j+1$ will be odd, and we will need to flip the array upside down. After the flip, the smallest allowable part is now in row $i-j+1$. Again, shifting $z \mapsto zq^{j+1}$ then makes the smallest allowable part equal to $j+2$. Thus, if $j$ does occur in the $(i-j+1)$-th row, then, in both cases, that adds a total of $zq^j P_{i-j+1}\bigl(zq^{j+1},q\bigr)$ to the generating function.

Now, if none of the parts $j$ in the $(2\ell-i+j)$-th rows actually occur, then we have a partition that could be counted by $P_{i+1}(zq,q)$. Again, to see this, we use the flipping trick: the smallest allowable part (the part 2) is in row $2\ell-i$, which is an even-indexed row (since $i$ is even). Flipping moves this part to row $i+1$. The partitions counted by $P_{i+1}(z,q)$ only are allowed to have 1 in row $i+1$; shifting $z\mapsto zq$ increases all parts by 1.

Since exactly one of the above cases ($j$ occurring in row $2\ell-i+j$ for $1 \le j \le i$, or none of those $j$s occurs), we obtain the functional equation $P_i(z,q) = \sum_{j=1}^i zq^j P_{i-j+1}\bigl(zq^{j+1},q\bigr)+ P_{i+1}(zq,q)$.

Finally, consider the case $i=\ell$. We follow one of the previous two cases, based on whether~$\ell$ is odd or even. The only difference is that if~$\ell$ is odd and none of the parts~$j$ located in rows $\ell-j+1$ (for $1\le j \le \ell$) occur, or if~$\ell$ is even and none of the parts~$j$ located in rows $\ell+j$ (for $1\le j \le \ell$) occur, then we have a partition counted by $P_{\ell}(zq,q)$ (as can be seen after flipping). This demonstrates $P_\ell(z,q) = \sum_{j=1}^\ell zq^j P_{\ell-j+1}\bigl(zq^{j+1},q\bigr) + P_\ell(zq,q)$.
\end{proof}

It will be more useful for us to rewrite these functional equations in the following form.
\begin{Theorem} \label{thm:streamlinedPeqs}
The generating functions $P_i(z,q)$ satisfy the following functional equations:
\begin{gather}P_0(z,q)-P_1(zq,q) =0,
\label{eq:newPrec1} \\
P_i(z,q) - P_{i+1}(zq,q) - P_{i-1}(zq,q) + (1-zq)P_{i}\bigl(zq^2,q\bigr) =0, \qquad 1 \le i \le \ell-1, \label{eq:newPrec2} \\
P_\ell(z,q) - P_\ell(zq,q) -P_{\ell-1}(zq,q) + (1-zq) P_{\ell}\bigl(zq^2,q\bigr) = 0. \label{eq:newPrec3}
\end{gather}
\end{Theorem}

\begin{proof}
Let $1\le i \le \ell$. We now take the $(i-1)$ instance of~\eqref{eq:origFE_gf1}, dilate $z\mapsto zq$, and then reindex $j$:
\begin{gather}
P_{i-1}(z,q) - P_{i}(zq,q) = \sum_{j=1}^{i-1} zq^j P_{i-j}\bigl(zq^{j+1},q\bigr), \notag \\
P_{i-1}(zq,q) - P_{i}\bigl(zq^2,q\bigr) = \sum_{j=1}^{i-1} zq^{j+1} P_{i-j}\bigl(zq^{j+2},q\bigr), \notag \\
P_{i-1}(zq,q) - P_{i}\bigl(zq^2,q\bigr) = \sum_{j=2}^{i} zq^{j} P_{i-j+1}\bigl(zq^{j+1},q\bigr), \qquad 1 \le i \le \ell. \label{eq:midFE3}
\end{gather}
Subtracting this from the $i$ instance of~\eqref{eq:origFE_gf1} gives
\begin{gather*}
P_i(z,q) - P_{i+1}(zq,q) - P_{i-1}(zq,q) + P_{i}\bigl(zq^2,q\bigr) \\
\qquad = \sum_{j=1}^i zq^j P_{i-j+1}\bigl(zq^{j+1},q\bigr) -\sum_{j=2}^{i} zq^{j} P_{i-j+1}\bigl(zq^{j+1},q\bigr) = zqP_i\bigl(zq^2,q\bigr), \label{eq:a7}
\end{gather*}
or
\begin{equation*}
P_i(z,q) - P_{i+1}(zq,q) - P_{i-1}(zq,q) + (1-zq)P_{i}\bigl(zq^2,q\bigr) =0, \qquad 1 \le i \le \ell-1.
\end{equation*}
Now, combining the $i=\ell$ instance of \eqref{eq:midFE3} with \eqref{eq:origFE_gf2} gives, similarly,
\begin{equation}
P_\ell(z,q) - P_\ell(zq,q) -P_{\ell-1}(zq,q) + (1-zq) P_{\ell}\bigl(zq^2,q\bigr) = 0. \tag*{\qed}
\end{equation}
\renewcommand{\qed}{}
\end{proof}

\section{Completing the proof: Andrews--Gordon case}\label{sec:comp}

Let us now define $T_i(z,q)$ ($0\le i \le \ell$) to be the right side of~\eqref{eq:mainresult1}:

\begin{equation} \label{eq:Tdef}
T_i(z,q)=\sum_{n_1,n_2,\dots,n_{\ell}\ge 0} \frac{z^{N_1}q^{N_1^2+N_2^2+\cdots+N_{\ell}^2+N_{i+1}+N_{i+2}+\cdots+N_{\ell}}}{(q;q)_{n_1}(q;q)_{n_2}\cdots(q;q)_{n_{\ell}}}.
\end{equation}

Our goal is to use the results of the previous section to show that $T_i(z,q)$ satisfies the following system of functional equations, corresponding to~\eqref{eq:newPrec1},~\eqref{eq:newPrec2}, and~\eqref{eq:newPrec3}:
\begin{gather}
T_0(z,q)-T_1(zq,q) =0, \label{eq:newTrec1} \\
T_i(z,q) - T_{i+1}(zq,q) - T_{i-1}(zq,q) + (1-zq)T_{i}\bigl(zq^2,q\bigr) =0, \qquad 1 \le i \le \ell-1, \label{eq:newTrec2} \\
T_\ell(z,q) - T_\ell(zq,q) -T_{\ell-1}(zq,q) + (1-zq) T_{\ell}\bigl(zq^2,q\bigr) = 0. \label{eq:newTrec3}
\end{gather}

As observed by Ole Warnaar,\footnote{Personal communication.} this set of functional equations exactly matches the (normalized) Corteel--Welsh~\cite{CW} functional equations for two-rowed cylindric partitions with profile $(c_1,c_2)$, where $c_1+c_2$ is odd (see their equation (3.5)). As such, we can complete the proof of Theorem~\ref{thm:mainresult1} by simply observing that we can recover~\eqref{eq:Tdef} by taking the limit as $L\to\infty$ in equations (7.24) and (7.25) of Warnaar's work~\cite{WarnCyl}. For the sake of making our paper self-contained (and because the proofs of Warnaar's formulas were omitted from~\cite{WarnCyl}), we will provide our own proof that $T_i(z,q)$ satisfies~\eqref{eq:newTrec1},~\eqref{eq:newTrec2}, and~\eqref{eq:newTrec3}.

Our proof involves what we will call the method of {\it atomic relations}, using computations similar to those of~\cite{BKRS,KanRus-cylindric}. The difference is that we will prove everything for general $\ell$ and $0\le i \le \ell$ without relying on computers to prove special cases. We note that the computations in this section were aided by Frank Garvan's {\tt qseries} Maple package~\cite{Gar}, along with the survey on Rogers--Ramanujan--Slater-type identities by James Mc Laughlin, Andrew V. Sills, and Peter Zimmer~\cite{MSZ}. To explain the concept of atomic relations, we first introduce some notation.

Let $\vv = \left\langle v_1,v_2,\dots,v_{\ell}\right\rangle$ and $\nn = \left\langle n_1,n_2,\dots,n_{\ell}\right\rangle$ be vectors in $\RR^\ell$. Then, define
\begin{equation*}
S_\vv(z,q)=\sum_{n_1,n_2,\dots,n_{\ell}\ge 0} \frac{z^{N_1}q^{N_1^2+N_2^2+\cdots+N_{\ell}^2+\vv \cdot \nn}}{(q;q)_{n_1}(q;q)_{n_2}\cdots(q;q)_{n_{\ell}}},
\end{equation*}
where $N_i = n_i +n_{i+1}+\cdots+n_{\ell}$ and $\vv\cdot \nn$ is the usual dot product. We will suppress the arguments and simply write $S_\vv$ for $S_\vv(z,q)$. To make our computations flow more smoothly, we make the following definitions:
\begin{gather*}
\zero = \langle 0,0,\dots,0\rangle, \\
\one = \langle 1,1,\dots,1\rangle, \\
\two = \langle 2,2,\dots,2\rangle, \\
\ee_i = \langle 0, \dots, 0, 1, 0, \dots, 0 \rangle, \qquad 1 \le i \le \ell, \\
\tr_i = \langle 0, \dots, 0, 1, 2, 3, \dots, \ell-i-1, \ell-i \rangle, \qquad 1 \le i \le \ell,
\end{gather*}
where, in $\ee_i$ and $\tr_i$, the 1 is located in position $i$. We will also define $\tr_{\ell+1} = \zero$.

Atomic relations are identities that state certain linear combinations of $S_\vv$ terms equal zero. These are generally obtained by multiplying the terms of the summation in $S_\vv$ by some factor that cancels out with one of factors in a $q$-Pochhammer symbol in the denominator, which then allows us to reindex the summation and then express the new summation again using $S_\vv$~notation. Our goal will be to rewrite the functional equations we would like to prove (in this case, \eqref{eq:newTrec1},~\eqref{eq:newTrec2}, and~\eqref{eq:newTrec3}) in terms of the $S_\vv$ notation, and then show that those new equations can themselves be written as linear combinations of our atomic relations. There is no guarantee that this method of atomic relations will work (in~\cite{KanRus-cylindric}, some negative examples of the method failing to prove certain conjectures are given), but when it does work (as in the cases in this present paper), it provides an organized way of keeping track of what may otherwise be rather complicated series manipulations. So now we have the following.
\begin{Lemma} \label{lem:atomAG}
For all $\vv$ and $1\le i \le \ell,$
\begin{equation*}
S_\vv - S_{\vv+\ee_i} = zq^{\vv\cdot\ee_i+i} S_{\vv+2\tr_1-2\tr_{i+1}}.
\end{equation*}
\end{Lemma}
\begin{proof}
We first multiply $S_\vv$ by $\bigl(1-q^{n_i}\bigr)$, and observe then that we have
\begin{equation*}
S_\vv-S_{\vv+\ee_i} =\sum_{n_1,n_2,\dots,n_{\ell}\ge 0} \frac{z^{N_1}q^{N_1^2+N_2^2+\cdots+N_{\ell}^2+\vv \cdot \nn}\bigl(1-q^{n_i}\bigr)}{(q;q)_{n_1}(q;q)_{n_2}\cdots(q;q)_{n_{\ell}}}.
\end{equation*}
Cancelling out a factor of $\bigl(1-q^{n_i}\bigr)$ in the denominator and reindexing $n_{i}-1 = \nihat$ produces
\begin{align*}
S_\vv-S_{\vv+\ee_i} &=\sum_{n_1,n_2,\dots,n_{\ell}\ge 0} \frac{z^{n_1+n_2+\cdots+n_{\ell}}q^{N_1^2+N_2^2+\cdots+N_{\ell}^2+\vv \cdot \nn}}{(q;q)_{n_1}(q;q)_{n_2}\cdots(q;q)_{n_{i}-1}\cdots(q;q)_{n_{\ell}}} \\
&=\sum_{n_1,n_2,\dots,\nihat,\dots,n_{\ell}\ge 0}\!\!\!\! \frac{z^{n_1+n_2+\cdots+\nihat+\cdots+n_{\ell}+1}q^{(N_1+1)^2+\cdots+(N_i+1)^2+N_{i+1}^2+\cdots+N_{\ell}^2+\vv \cdot \nn+\vv\cdot\ee_i}}{(q;q)_{n_1}(q;q)_{n_2}\cdots(q;q)_{\nihat}\cdots(q;q)_{n_{\ell}}} \\
&=zq^{\vv\cdot\ee_i+i}\sum_{n_1,n_2,\dots,\nihat,\dots,n_{\ell}\ge 0} \frac{z^{n_1+n_2+\cdots+\nihat+\cdots+n_{\ell}}q^{N_1^2+\cdots+N_{\ell}^2+2(N_1+\cdots+N_i)+\vv \cdot \nn}}{(q;q)_{n_1}(q;q)_{n_2}\cdots(q;q)_{\nihat}\cdots(q;q)_{n_{\ell}}} \\
&= zq^{\vv\cdot\ee_i+i} S_{\vv+2\tr_1-2\tr_{i+1}}.
\end{align*}
(In a slight abuse of notation, we used the notation $N_j$ and $\nn$ with $\nihat$ replacing $n_i$.)
\end{proof}

So we now define our atomic relations: for $1\le i \le \ell$,
\begin{equation*}
\rel^i_\vv:=S_\vv - S_{\vv+\ee_i} - zq^{\vv\cdot\ee_i+i} S_{\vv+2\tr_1-2\tr_{i+1}}.
\end{equation*}
Each $\rel^i_\vv$ is, of course, equal to zero. If we can obtain a certain expression as a linear combination of $\rel^i_\vv$s, then that expression will also be equal to zero.

We now translate the desired identities~\eqref{eq:newTrec1}--\eqref{eq:newTrec3} into $S_\vv$ notation. Recalling~\eqref{eq:Tdef}, we see the identification $T_i(z,q) = S_{\tr_{i+1}}(z,q).$ We obtain
\begin{gather*}
S_{\tr_1}(z,q)-S_{\tr_2}(zq,q) =0, \\
S_{\tr_{i+1}}(z,q) - S_{\tr_{i+2}}(zq,q) - S_{\tr_{i}}(zq,q) + (1-zq)S_{\tr_{i+1}}\bigl(zq^2,q\bigr) =0, \qquad 1 \le i \le \ell-1, \\
S_{\zero}(z,q) - S_{\zero}(zq,q) -S_{\ee_\ell}(zq,q) + (1-zq) S_{\zero}\bigl(zq^2,q\bigr) = 0.
\end{gather*}
Note that $S_\vv(zq,q) = S_{\vv+\one} (z,q)$ and $S_\vv\bigl(zq^2,q\bigr) = S_{\vv+\two} (z,q)$. This allows us to rewrite these as
\begin{gather}
S_{\tr_1}(z,q)-S_{\one+\tr_2}(z,q) =0, \label{eq:finFE1} \\
S_{\tr_{i+1}}(z,q) - S_{\one+\tr_{i+2}}(z,q) - S_{\one+\tr_{i}}(z,q) + (1-zq)S_{\two+\tr_{i+1}}(z,q) =0, \qquad\! 1 \le i \le \ell-1,\!\!\! \label{eq:finFE2} \\
S_{\zero}(z,q) - S_{\one}(z,q) -S_{\one+\ee_\ell}(z,q) + (1-zq) S_{\two}(z,q) = 0. \label{eq:finFE3}
\end{gather}
So our goal is to obtain~\eqref{eq:finFE1},~\eqref{eq:finFE2}, and~\eqref{eq:finFE3} as linear combinations of atomic relations. To that end, we need two propositions. (Appendix~\ref{subsec:sercomp} explicitly shows how to obtain Propositions~\ref{prop:comp1} and~\ref{prop:comp2} in the case $\ell=4$.)
\begin{Proposition} \label{prop:comp1}
$S_{\zero} - S_{\one} -
S_{\one+\ee_\ell} + (1-zq) S_{\two}=0.$
\end{Proposition}
\begin{proof}
We can obtain this as the following linear combinations of relations:
\begin{align*}
&\sum_{i=1}^\ell \rel^i_{\ee_1+\cdots+\ee_{i-1}} - \sum_{i=1}^{\ell-1}\rel^i_{\two-\ee_1-\cdots-\ee_i} \\
&=\sum_{i=1}^\ell \bigl(
S_{\ee_1+\cdots+\ee_{i-1}} - S_{\ee_1+\cdots+\ee_i} - zq^{i} S_{\ee_1+\cdots+\ee_{i-1}+2\tr_1-2\tr_{i+1}}\bigr) \\
&\quad- \sum_{i=1}^{\ell-1}\bigl(
S_{\two-\ee_1-\cdots-\ee_i} - S_{\two-\ee_1-\cdots-\ee_{i-1}} - zq^{i+1} S_{\two-\ee_1-\cdots-\ee_i+2\tr_1-2\tr_{i+1}}\bigr) \\
&=S_{\zero} - S_{\one} - \sum_{i=1}^\ell zq^{i} S_{\ee_1+\cdots+\ee_{i-1}+2\tr_1-2\tr_{i+1}} -
S_{\one+\ee_\ell} + S_{\two} + \sum_{i=1}^{\ell-1}zq^{i+1} S_{\two-\ee_1-\cdots-\ee_i+2\tr_1-2\tr_{i+1}} \\
&= S_{\zero} - S_{\one} -
S_{\one+\ee_\ell} + S_{\two} - \sum_{i=1}^\ell zq^{i} S_{\ee_1+\cdots+\ee_{i-1}+2\tr_1-2\tr_{i+1}} + \sum_{i=2}^{\ell}zq^{i} S_{\two-\ee_1-\cdots-\ee_{i-1}+2\tr_1-2\tr_{i}} \\
&= S_{\zero} - S_{\one} -
S_{\one+\ee_\ell} + S_{\two} -
 zq^{1} S_{2\tr_1-2\tr_{2}} \\
&\quad-\sum_{i=2}^\ell zq^{i} \bigl(S_{\ee_1+\cdots+\ee_{i-1}+2\tr_1-2\tr_{i+1}} -S_{\two-\ee_1-\cdots-\ee_{i-1}+2\tr_1-2\tr_{i}}\bigr) \\
&= S_{\zero} - S_{\one} -
S_{\one+\ee_\ell} + S_{\two} -
 zq^{1} S_{\two} -0 \\
&= S_{\zero} - S_{\one} - S_{\one+\ee_\ell} + (1-zq)S_{\two} \\
&= 0.
\end{align*}
Above, we used the equality
\begin{equation*}
\ee_1+\cdots+\ee_{i-1}+2\tr_1-2\tr_{i+1}={\two-\ee_1-\cdots-\ee_{i-1}+2\tr_1-2\tr_{i}}
\end{equation*} which follows from $
2\bigl(\ee_1+\cdots+\ee_{i-1}\bigr)=\two-2\tr_{i}+2\tr_{i+1}.$
\end{proof}

\begin{Proposition} \label{prop:comp2}
For $1 \le j \le \ell-1$,
\begin{equation*}
S_{\tr_{j+1}} -S_{\ee_1+\cdots+\ee_{j}+\tr_{j+1}}- S_{\two-\ee_1-\cdots-\ee_{j-1}+\tr_{j+1}} + (1-zq)S_{\two+\tr_{j+1}}=0.
\end{equation*}
\end{Proposition}
\begin{proof}
\begin{align*}
&\sum_{i=1}^j \rel^i_{\ee_1+\cdots+\ee_{i-1}+\tr_{j+1}} - \sum_{i=1}^{j-1}\rel^i_{\two-\ee_1-\cdots-\ee_i+\tr_{j+1}} \\
&=\sum_{i=1}^j \bigl( S_{\ee_1+\cdots+\ee_{i-1}+\tr_{j+1}} - S_{\ee_1+\cdots+\ee_{i}+\tr_{j+1}} - zq^{i} S_{\ee_1+\cdots+\ee_{i-1}+\tr_{j+1}+2\tr_1-2\tr_{i+1}}\bigr) \\
& \quad - \sum_{i=1}^{j-1}\bigl( S_{\two-\ee_1-\cdots-\ee_i+\tr_{j+1}} - S_{\two-\ee_1-\cdots-\ee_{i-1}+\tr_{j+1}} - zq^{i+1} S_{\two-\ee_1-\cdots-\ee_i+\tr_{j+1}+2\tr_1-2\tr_{i+1}}\bigr) \\
&=S_{\tr_{j+1}} -S_{\ee_1+\cdots+\ee_{j}+\tr_{j+1}} - \sum_{i=1}^j zq^{i} S_{\ee_1+\cdots+\ee_{i-1}+\tr_{j+1}+2\tr_1-2\tr_{i+1}} \\
& \quad - S_{\two-\ee_1-\cdots-\ee_{j-1}+\tr_{j+1}} + S_{\two+\tr_{j+1}} +\sum_{i=1}^{j-1}zq^{i+1} S_{\two-\ee_1-\cdots-\ee_i+\tr_{j+1}+2\tr_1-2\tr_{i+1}}
\\
&=S_{\tr_{j+1}} -S_{\ee_1+\cdots+\ee_{j}+\tr_{j+1}}- S_{\two-\ee_1-\cdots-\ee_{j-1}+\tr_{j+1}} + S_{\two+\tr_{j+1}} \\
& \quad - \sum_{i=1}^j zq^{i} S_{\ee_1+\cdots+\ee_{i-1}+\tr_{j+1}+2\tr_1-2\tr_{i+1}} +\sum_{i=2}^{j}zq^{i} S_{\two-\ee_1-\cdots-\ee_{i-1}+\tr_{j+1}+2\tr_1-2\tr_{i}} \\
&=S_{\tr_{j+1}} -S_{\ee_1+\cdots+\ee_{j}+\tr_{j+1}}- S_{\two-\ee_1-\cdots-\ee_{j-1}+\tr_{j+1}} + S_{\two+\tr_{j+1}} - zq^{1} S_{\tr_{j+1}+2\tr_1-2\tr_{2}} \\
& \quad - \sum_{i=1}^j zq^{i} \bigl( S_{\ee_1+\cdots+\ee_{i-1}+\tr_{j+1}+2\tr_1-2\tr_{i+1}} - S_{\two-\ee_1-\cdots-\ee_{i-1}+\tr_{j+1}+2\tr_1-2\tr_{i}}\bigr) \\
&=S_{\tr_{j+1}} -S_{\ee_1+\cdots+\ee_{j}+\tr_{j+1}}- S_{\two-\ee_1-\cdots-\ee_{j-1}+\tr_{j+1}} + (1-zq)S_{\two+\tr_{j+1}} \\
&=0.
\end{align*}
Again, we used the equality
\begin{equation*}
\ee_1+\cdots+\ee_{i-1}+\tr_{j+1}+2\tr_1-2\tr_{i+1}=
\two-\ee_1-\cdots-\ee_{i-1}+\tr_{j+1}+2\tr_1-2\tr_{i}
\end{equation*}
which follows from $2(\ee_1+\cdots+\ee_{i-1})-2\tr_{i+1} =\two-2\tr_{i}.$ We also used the equality $\two+\tr_{j+1} ={\tr_{j+1}+2\tr_1-2\tr_{2}}$.
\end{proof}

So now Proposition~\ref{prop:comp1} gives us~\eqref{eq:finFE3} verbatim. Equation~\eqref{eq:finFE1} follows from the fact that $\tr_1=\one+\tr_2$. Also, we can use Proposition~\ref{prop:comp2}, together with the facts that $\ee_1+\cdots+\ee_{i}+\tr_{i+1} = \one+\tr_{i+2}$ and $\one+\tr_{i} =\two-\ee_1-\cdots-\ee_{i-1}+\tr_{i+1}$, to obtain~\eqref{eq:finFE2}. Finally, checking the initial conditions ($P_i(0,q)=T_i(0,q)=1$) allows us to complete our proof of Theorem~\ref{thm:mainresult1}.

\section{Partitions: Andrews--Bressoud case}\label{sec:partitions4}

Fix $\ell\ge1$. In this section (and the next), we construct two different types of (extended) frequency arrays on $2\ell-1$ rows.

What we will call an ``odd-odd'' frequency array will have odd numbers in the first and last rows. Its second column will consist of all $\cdot$'s (again representing parts that are forbidden from appearing). The entries of the first column in rows (with odd index) 1 through $\ell$ will be $k_i$, while the remainder of the entries in the first column (in rows with odd indices $2\left\lceil \ell/2\right\rceil+1$ through $2\ell-1$) will have $\cdot$s in them. Here is the general form of these:
\begin{equation*}
\begin{matrix}
k_1 && 1 && 3 && 5 && 7 && 9 && 11 & \\
& \cdot && 2 && 4 && 6 && 8 && 10 && 12\\
k_3 && 1 && 3 && 5 && 7 && 9 && 11 & \\
& \cdot && 2 && 4 && 6 && 8 && 10 && 12\\
k_5 && 1 && 3 && 5 && 7 && 9 && 11 & \\
& \cdot && 2 && 4 && 6 && 8 && 10 && 12\\
\vdots & \vdots &\vdots& \vdots &\vdots& \vdots &\vdots& \vdots &\vdots& \vdots &\vdots& \vdots\\
k_{2\left\lceil \ell/2\right\rceil-1} && 1 && 3 && 5 && 7 && 9 && 11 & \\
& \cdot && 2 && 4 && 6 && 8 && 10 && 12\\
\cdot && 1 && 3 && 5 && 7 && 9 && 11 & \\
& \cdot && 2 && 4 && 6 && 8 && 10 && 12\\
\vdots & \vdots &\vdots& \vdots &\vdots& \vdots &\vdots& \vdots &\vdots& \vdots &\vdots& \vdots\\
\cdot && 1 && 3 && 5 && 7 && 9 && 11 & \\
& \cdot && 2 && 4 && 6 && 8 && 10 && 12 \\
\cdot && 1 && 3 && 5 && 7 && 9 && 11 &
 \end{matrix}\quad\dots \,.
\end{equation*}
Note that $2\left\lceil \ell/2\right\rceil-1$ is the largest odd number that is at most $\ell$.

Correspondingly, what we will call an ``even-even'' frequency array will have even numbers in the first and last rows. Its second column will have $k_0$ in the first row, and then the rest of the entries will be $\cdot$s. The entries of the first column in rows (with even index) 2 through $\ell$ will be~$k_i$, while the remainder of the entries in the first column (in rows with even indices $2\left\lfloor \ell/2\right\rfloor+2$ through $2\ell-2$) will have $\cdot$s in them. Again, we show the general form of these:
\begin{equation*}
\begin{matrix}
& k_0 && 2 && 4 && 6 && 8 && 10 && 12\\
k_2 && 1 && 3 && 5 && 7 && 9 && 11 & \\
& \cdot && 2 && 4 && 6 && 8 && 10 && 12\\
k_4 && 1 && 3 && 5 && 7 && 9 && 11 & \\
& \cdot && 2 && 4 && 6 && 8 && 10 && 12\\
\vdots & \vdots &\vdots& \vdots &\vdots& \vdots &\vdots& \vdots &\vdots& \vdots &\vdots& \vdots\\
k_{2\left\lfloor \ell/2\right\rfloor} && 1 && 3 && 5 && 7 && 9 && 11 & \\
& \cdot && 2 && 4 && 6 && 8 && 10 && 12\\
\cdot && 1 && 3 && 5 && 7 && 9 && 11 & \\
& \cdot && 2 && 4 && 6 && 8 && 10 && 12\\
\vdots & \vdots &\vdots& \vdots &\vdots& \vdots &\vdots& \vdots &\vdots& \vdots &\vdots& \vdots\\
\cdot && 1 && 3 && 5 && 7 && 9 && 11 & \\
& \cdot && 2 && 4 && 6 && 8 && 10 && 12 \\
\cdot && 1 && 3 && 5 && 7 && 9 && 11 & \\
& \cdot && 2 && 4 && 6 && 8 && 10 && 12
 \end{matrix}\quad\dots \,.
\end{equation*}
Note that $2\left\lfloor \ell/2\right\rfloor$ is the largest even number that is at most $\ell$. Here is an example of an odd-odd frequency array for $\ell=6$:
\begin{equation*}
\begin{matrix}
k_1 && f_1 && f_3 && f_5 && f_7 && f_9 && f_{11} & \\
& \cdot && f_2 && f_4 && f_6 && f_8 && f_{10} && f_{12}\\
k_3 && f_1 && f_3 && f_5 && f_7 && f_9 && f_{11} & \\
& \cdot && f_2 && f_4 && f_6 && f_8 && f_{10} && f_{12}\\
k_5 && f_1 && f_3 && f_5 && f_7 && f_9 && f_{11} & \\
& \cdot && f_2 && f_4 && f_6 && f_8 && f_{10} && f_{12}\\
\cdot && f_1 && f_3 && f_5 && f_7 && f_9 && f_{11} & & \\
& \cdot && f_2 && f_4 && f_6 && f_8 && f_{10} && f_{12}\\
\cdot && f_1 && f_3 && f_5 && f_7 && f_9 && f_{11} & \\
& \cdot && f_2 && f_4 && f_6 && f_8 && f_{10} && f_{12}\\
\cdot && f_1 && f_3 && f_5 && f_7 && f_9 && f_{11} & \\
\end{matrix}\quad\dots \,.
\end{equation*}
And here is an example of an even-even frequency array for $\ell=6$:
\begin{equation*}
\begin{matrix}
& k_0 && f_2 && f_4 && f_6 && f_8 && f_{10} && f_{12}\\
k_2 && f_1 && f_3 && f_5 && f_7 && f_9 && f_{11} & \\
& \cdot && f_2 && f_4 && f_6 && f_8 && f_{10} && f_{12}\\
k_4 && f_1 && f_3 && f_5 && f_7 && f_9 && f_{11} & \\
& \cdot && f_2 && f_4 && f_6 && f_8 && f_{10} && f_{12}\\
k_6 && f_1 && f_3 && f_5 && f_7 && f_9 && f_{11} & \\
& \cdot && f_2 && f_4 && f_6 && f_8 && f_{10} && f_{12}\\
\cdot && f_1 && f_3 && f_5 && f_7 && f_9 && f_{11} & \\
& \cdot && f_2 && f_4 && f_6 && f_8 && 10 && f_{12}\\
\cdot && f_1 && f_3 && f_5 && f_7 && f_9 && f_{11} & \\
& \cdot && f_2 && f_4 && f_6 && f_8 && f_{10} && f_{12}\\
\end{matrix}\quad\dots \,.
\end{equation*}

For odd $1 \le i \le \ell$, define $P^\star_i(z,q)$ to be the bivariate generating function for $[k_0,k_1\dots,k_\ell]=[0,\dots,0,1,0,\dots,0]$-admissible partitions (using odd-odd frequency arrays) of $n$, and for even ${0 \le i \le \ell}$, define $P^\star_i(z,q)$ to be the bivariate generating function for $[k_0,k_1\dots,k_\ell]=[0,\dots,0,1,\allowbreak 0,\dots,0]$-admissible partitions (using even-even frequency arrays) of~$n$. As before, we define admissible partitions through downward paths: since $k=1$, any downward path in these frequency arrays must have at most one part. In both cases, the 1 in $[0,\dots,0,1,0,\dots,0]$ is located in position $i$. We will now show that the bivariate generating functions for the admissible partitions related to these arrays satisfy the following system of functional equations.
\begin{Theorem}\label{thm:funceq1B}
The generating functions $P^\star_i(z,q)$ satisfy the following functional equations:
\begin{gather*}
P^\star_i(z,q) - P^\star_{i+1}(zq,q) = \sum_{j=1}^i zq^j P^\star_{i-j+1}\bigl(zq^{j+1},q\bigr), \qquad 0 \le i \le \ell-1, \\
P^\star_\ell(z,q) - P^\star_{\ell-1}(zq,q) = \sum_{j=1}^\ell zq^j P^\star_{\ell-j+1}\bigl(zq^{j+1},q\bigr).
\end{gather*}
\end{Theorem}

\begin{proof} Unlike in the proof of Theorem~\ref{thm:funceq1A}, the ``flipping'' trick is not required in any of these cases.

First, consider the case $i=0$. We present the extended frequency arrays for $P^\star_0(z,q)$ and $P^\star_1(z,q)$:
\begin{gather*}
P^\star_0(z,q)\colon \quad \begin{matrix}
& \ficone && 2 && 4 && 6 && 8 && 10 && 12\\
\cdot && \cdot && 3 && 5 && 7 && 9 && 11 & \\
& \cdot && \cdot && 4 && 6 && 8 && 10 && 12\\
\cdot && \cdot && \cdot && 5 && 7 && 9 && 11 & \\
& \cdot && \cdot && \cdot && 6 && 8 && 10 && 12\\
\vdots & \vdots &\vdots& \vdots &\vdots& \vdots &\vdots& \vdots &\vdots& \vdots &\vdots& \vdots&\vdots& \vdots\vspace{3mm}
 \end{matrix}\quad\dots\,, \\
P^\star_1(z,q)\colon \quad
\begin{matrix}
\ficone && 1 && 3 && 5 && 7 && 9 && 11 & \\
& \cdot && 2 && 4 && 6 && 8 && 10 && 12\\
\cdot && \cdot && 3 && 5 && 7 && 9 && 11 & \\
& \cdot && \cdot && 4 && 6 && 8 && 10 && 12\\
\cdot && \cdot && \cdot && 5 && 7 && 9 && 11 & \\
& \cdot && \cdot && \cdot && 6 && 8 && 10 && 12\\
\vdots & \vdots &\vdots& \vdots &\vdots& \vdots &\vdots& \vdots &\vdots& \vdots &\vdots& \vdots&\vdots& \vdots\\
\end{matrix}\quad\dots\,.
\end{gather*}
The generating function $P^\star_1(zq,q)$ corresponds to adding~1 to every part in the partitions counted by $P^\star_1(z,q)$, and so it is easy to see that $P^\star_0(z,q)=P^\star_1(zq,q)$.

Next, consider the case where $i\not \in\{0,\ell\}$. Suppose that we have an admissible partition counted by $P^\star_i(z,q)$. Consider, for $1 \le j \le i$, the set of parts $j$ that may occur in row $i-j+1$; at most one of these parts can occur in this partition. Once again, let us look at a similar picture to the one we presented in the proof of Theorem~\ref{thm:funceq1A}:
\begin{gather*}
P_i^\star(z,q)\colon \quad
\begin{matrix}
\!\cdots\! && \!\cdots\! && \!\cdots\! && \cdot && \color{blue}{f_{i}} && f_{i+2} && f_{i+4} & \\
& \!\cdots\! && \!\cdots\! && \cdot && \color{blue}{f_{i-1}} && f_{i+1} && f_{i+3} && f_{i+5}\\
\vdots & \vdots &\vdots& \vdots &\vdots& \vdots &\iddots& \vdots &\vdots& \vdots &\vdots& \vdots&\vdots& \vdots\\
& \cdot && \cdot && \color{blue}{f_4} && f_6 && f_8 && f_{10} && f_{12}\\
\cdot &&\cdot && \color{blue}{f_3} && f_5 && f_7 && f_9 && f_{11} & \\
& \cdot && \color{blue}{f_2} && f_4 && f_6 && f_8 && f_{10} && f_{12}\\
\ficone &&\color{blue}{f_1} && f_3 && f_5 && f_7 && f_9 && f_{11} & \\
& \cdot && f_2 && f_4 && f_6 && f_8 && f_{10} && f_{12}\\
\cdot &&\cdot && f_3 && f_5 && f_7 && f_9 && f_{11} & \\
& \cdot && \cdot && f_4 && f_6 && f_8 && f_{10} && f_{12}\\
\vdots & \vdots &\vdots& \vdots &\vdots& \vdots &\vdots& \vdots &\vdots& \vdots &\vdots& \vdots&\vdots& \vdots
 \end{matrix}\quad\dots\,.
\end{gather*}
If $i$ is odd, then this will be an odd-odd frequency array (with odd numbers in the first and last rows), while if $i$ is even, this will be an even-even frequency array. Suppose that $j$ does occur in the $(i-j+1)$-th row. After removing this part $j$, the remaining parts will be counted by \smash{$P^\star_{i-j+1}\bigl(zq^{j+1},q\bigr)$}. To see this, note that the smallest allowable part, $j+2$, is in row $i-j+1$. Recall that the only row in which the part 1 can occur in the frequency array for \smash{$P^\star_{i-j+1}(z,q)$} is row $i-j+1$, and so shifting $z\mapsto zq^{j+1}$ gives the desired frequency array. (Note that, if $j+1$ is odd, then shifting $z\mapsto zq^{j+1}$ switches an odd-odd frequency array to an even-even frequency array or vice versa.) Thus, if $j$ does occur in the $(i-j+1)$-th row, then that adds a total of \smash{$zq^j P^\star_{i-j+1}\bigl(zq^{j+1},q\bigr)$} to the generating function, where the factor of $zq^j$ comes from the part $j$.
Now, if none of the parts $j$ in the $(i-j+1)$-th rows actually occur, then we have a partition that could be counted by $P^\star_{i+1}(zq,q)$. The only row that partitions counted by $P^\star_{i+1}(z,q)$ are allowed to have the part 1 occur in is row $i+1$; shifting $z\mapsto zq$ increases all parts by 1. (Note that this shift switches an odd-odd frequency array to an even-even frequency array or vice versa.) Since exactly one of the above cases ($j$ occurring in row $i-j+1$ for $1 \le j \le i$, or none of those $j$s occurs), we obtain the functional equation $P^\star_i(z,q) = \sum_{j=1}^i zq^j P^\star_{i-j+1}\bigl(zq^{j+1},q\bigr)+ P^\star_{i+1}(zq,q)$.

Finally, consider the case $i=\ell$. We follow the previous case. The only difference is that if none of the parts $j$ in row $i-j+1$ occur, then we have a partition that is counted by \smash{$P^\star_{\ell-1}(zq,q)$}. This demonstrates \smash{$P^\star_\ell(z,q) = \sum_{j=1}^\ell zq^j P^\star_{\ell-j+1}\bigl(zq^{j+1},q\bigr) + P^\star_\ell(zq,q)$}.
\end{proof}

We illustrate the $\ell=4$ case of this theorem in Appendix~\ref{subsec:combspec2}. Once again, we can follow a~very similar process to that of the proof of Theorem~\ref{thm:streamlinedPeqs} to obtain the following.
\begin{Theorem}
The generating functions $P^\star_i(z,q)$ satisfy the following functional equations:
\begin{gather}
P^\star_0(z,q)-P^\star_1(zq,q) =0,\label{eq:newPsarrec1} \\
P^\star_i(z,q) - P^\star_{i+1}(zq,q) - P^\star_{i-1}(zq,q) + (1-zq)P^\star_{i}\bigl(zq^2,q\bigr) =0, \qquad 1 \le i \le \ell-1, \label{eq:newPsarrec2} \\
P^\star_\ell(z,q) - 2P^\star_{\ell-1}(zq,q) + (1-zq) P^\star_{\ell}\bigl(zq^2,q\bigr) = 0.\label{eq:newPsarrec3}
\end{gather}
\end{Theorem}
Notice that~\eqref{eq:newPsarrec1} and~\eqref{eq:newPsarrec2} are essentially the same as~\eqref{eq:newPrec1} and~\eqref{eq:newPrec2}; the only difference between this theorem and Theorem~\ref{thm:streamlinedPeqs} is the change from~\eqref{eq:newPrec3} to~\eqref{eq:newPsarrec3}.

\section{Completing the proof: Andrews--Bressoud case}\label{sec:comp5}

Define (for $0 \le i \le \ell$)
\begin{equation} \label{eq:Tstardefn}
T^\star_i(z,q)=\sum_{n_1,n_2,\dots,n_{\ell}\ge 0} \frac{z^{N_1}q^{N_1^2+N_2^2+\cdots+N_{\ell}^2+N_{i+1}+N_{i+2}+\cdots+N_{\ell}}}{(q;q)_{n_1}(q;q)_{n_2} \cdots(q;q)_{n_{\ell-1}}\bigl(q^2;q^2\bigr)_{n_{\ell}}}.
\end{equation}
As before, we want to show that $T^\star_i(z,q)$ satisfies the following set of functional equations:
\begin{gather}
T^\star_0(z,q)-T^\star_1(zq,q) =0, \label{eq:newTsarrec1} \\
T^\star_i(z,q) - T^\star_{i+1}(zq,q) - T^\star_{i-1}(zq,q) + (1-zq)T^\star_{i}\bigl(zq^2,q\bigr) =0, \qquad 1 \le i \le \ell-1, \label{eq:newTsarrec2} \\
T^\star_\ell(z,q) - 2T^\star_{\ell-1}(zq,q) + (1-zq) T^\star_{\ell}\bigl(zq^2,q\bigr) = 0. \label{eq:newTsarrec3}
\end{gather}
And once again, these are the same as the normalized Corteel--Welsh functional equations~\cite{CW} for two-rowed cylindric partitions -- this time, for profiles $(c_1,c_2)$ where $c_1+c_2$ is even. Again, we can simply appeal to Section~7 of Warnaar's work~\cite{WarnCyl}, but we choose to give the full proof. Analogously to before, we define
\begin{equation*}
S^\star_\vv(z,q)=\sum_{n_1,n_2,\dots,n_{\ell}\ge 0} \frac{z^{N_1}q^{N_1^2+N_2^2+\cdots+N_{\ell}^2+\vv \cdot \nn}}{(q;q)_{n_1}(q;q)_{n_2}\cdots\bigl(q^2;q^2\bigr)_{n_{\ell}}}.
\end{equation*}

\begin{Lemma}[atomic relations]
For all $\vv$ and $1\le i \le \ell-1,$
\begin{equation} \label{eq:a4}
S^\star_\vv - S^\star_{\vv+\ee_i} = zq^{\vv\cdot\ee_i+i} S^\star_{\vv+2\tr_1-2\tr_{i+1}}.
\end{equation}
Also,
\begin{equation} \label{eq:a5}
S^\star_\vv - S^\star_{\vv+2\ee_\ell} = zq^{\vv\cdot\ee_\ell+\ell} S^\star_{\vv+2\tr_1}.
\end{equation}
\end{Lemma}
\begin{proof}
As before: for $1\le i \le \ell-1$, multiply $S^\star_\vv$ by $\bigl(1-q^{n_i}\bigr)$, cancel out a factor in the denominator, and reindex with respect to $n_i$ to obtain~\eqref{eq:a4}. Also, multiply $S^\star_\vv$ by $\bigl(1-q^{2n_\ell}\bigr)$ and do the same to obtain~\eqref{eq:a5}.
\end{proof}

So now define, for $1\le i \le \ell-1$,
\begin{equation*}
\rel^i_\vv:=S^\star_\vv - S^\star_{\vv+\ee_i} - zq^{\vv\cdot\ee_i+i} S^\star_{\vv+2\tr_1-2\tr_{i+1}},
\end{equation*}
and for $i=\ell$,
\begin{equation*}
\rel^\ell_\vv:=S^\star_\vv - S^\star_{\vv+2\ee_\ell} - zq^{\vv\cdot\ee_\ell+\ell} S^\star_{\vv+2\tr_1}.
\end{equation*}
The atomic relations are nearly the same as those in Section~\ref{sec:comp}. The only difference is that $\rel^\ell_\vv$ has changed.

\begin{Proposition} \label{prop:comp2_AB}
For $1 \le j \le \ell-1$,
\begin{equation*}
S^\star_{\tr_{j+1}} -S^\star_{\ee_1+\cdots+\ee_{j}+\tr_{j+1}}- S^\star_{\two-\ee_1-\cdots-\ee_{j-1}+\tr_{j+1}} + (1-zq)S^\star_{\two+\tr_{j+1}}=0.
\end{equation*}
\end{Proposition}
\begin{proof}
This is the equivalent proposition to the previous Proposition~\ref{prop:comp2}, but with $S$ replaced by $S^\star$. Since that one only relied on $\rel^i_\vv$ for $1\le i \le \ell-1$, which are exactly the same in both cases, the proof is exactly the same.
\end{proof}

\begin{Proposition} \label{prop:comp1_AB}
We have $S^\star_{\zero} - 2S^\star_{\one+\ee_\ell} + (1-zq) S^\star_{\two}=0$.
\end{Proposition}
\begin{proof}
We can obtain this as the following linear combinations of relations:
\begin{align*}
&\sum_{i=1}^\ell \rel^i_{\ee_1+\cdots+\ee_{i-1}} - \sum_{i=1}^{\ell-1}\rel^i_{\two-\ee_1-\cdots-\ee_i} \\
&=\sum_{i=1}^{\ell-1} \bigl(
S^\star_{\ee_1+\cdots+\ee_{i-1}} - S^\star_{\ee_1+\cdots+\ee_i} - zq^{i} S^\star_{\ee_1+\cdots+\ee_{i-1}+2\tr_1-2\tr_{i+1}}\bigr) \\
&\quad+S^\star_{\ee_1+\cdots+\ee_{\ell-1}} - S^\star_{{\ee_1+\cdots+\ee_{\ell-1}}+2\ee_\ell} - zq^{\ell} S^\star_{{\ee_1+\cdots+\ee_{\ell-1}}+2\tr_1} \\
&\quad- \sum_{i=1}^{\ell-1}\bigl(
S^\star_{\two-\ee_1-\cdots-\ee_i} - S^\star_{\two-\ee_1-\cdots-\ee_{i-1}} - zq^{i+1} S^\star_{\two-\ee_1-\cdots-\ee_i+2\tr_1-2\tr_{i+1}}\bigr) \\
&=
S^\star_{\zero} - S^\star_{\one+\ee_\ell} - \sum_{i=1}^{\ell} zq^{i} S^\star_{\ee_1+\cdots+\ee_{i-1}+2\tr_1-2\tr_{i+1}} -
S^\star_{\one+\ee_\ell}\\
&\quad + S^\star_{\two} + \sum_{i=1}^{\ell-1}zq^{i+1} S^\star_{\two-\ee_1-\cdots-\ee_i+2\tr_1-2\tr_{i+1}} \\
&=
S^\star_{\zero} -2S^\star_{\one+\ee_\ell} + S^\star_{\two} - \sum_{i=1}^\ell zq^{i} S^\star_{\ee_1+\cdots+\ee_{i-1}+2\tr_1-2\tr_{i+1}} + \sum_{i=2}^{\ell}zq^{i} S^\star_{\two-\ee_1-\cdots-\ee_{i-1}+2\tr_1-2\tr_{i}} \\
&= S^\star_{\zero} -2S^\star_{\one+\ee_\ell} + S^\star_{\two} -
 zq^{1} S^\star_{2\tr_1-2\tr_{2}}\\
 &\quad -\sum_{i=2}^\ell zq^{i} \bigl(S^\star_{\ee_1+\cdots+\ee_{i-1}+2\tr_1-2\tr_{i+1}} -S^\star_{\two-\ee_1-\cdots-\ee_{i-1}+2\tr_1-2\tr_{i}}\bigr) \\
&= S^\star_{\zero} - 2S^\star_{\one+\ee_\ell} + S^\star_{\two} -
 zq^{1} S^\star_{\two} \\
&= S^\star_{\zero} - 2S^\star_{\one+\ee_\ell} + (1-zq)S^\star_{\two} \\
&=0.
\end{align*}
Once again, we used the facts that \[{\ee_1+\cdots+\ee_{i-1}+2\tr_1-2\tr_{i+1}}={\two-\ee_1-\cdots-\ee_{i-1}+2\tr_1-2\tr_{i}}\] and $2\tr_1-2\tr_2=\two$.
\end{proof}

Let us translate~\eqref{eq:newTsarrec1}--\eqref{eq:newTsarrec3} into $S^\star_\vv$ notation. Recalling~\eqref{eq:Tstardefn}, we see the identification $T^\star_i(z,q) = S^\star_{\tr_{i+1}}(z,q)$, producing
\begin{gather*}
S^\star_{\tr_1}(z,q)-S^\star_{\tr_2}(zq,q) =0, \\
S^\star_{\tr_{i+1}}(z,q) - S^\star_{\tr_{i+2}}(zq,q) - S^\star_{\tr_{i}}(zq,q) + (1-zq)S^\star_{\tr_{i+1}}\bigl(zq^2,q\bigr) =0, \qquad 1 \le i \le \ell-1, \\
S^\star_{\zero}(z,q) -2S^\star_{\ee_\ell}(zq,q) + (1-zq) S^\star_{\zero}\bigl(zq^2,q\bigr) = 0.
\end{gather*}

Note that $S^\star_\vv(zq,q) = S^\star_{\vv+\one} (z,q)$ and $S^\star_\vv\bigl(zq^2,q\bigr) = S^\star_{\vv+\two} (z,q)$. This allows us to rewrite these~as
\begin{gather}
S^\star_{\tr_1}(z,q)-S^\star_{\one+\tr_2}(z,q) =0,\label{eq:finFE_AB1} \\
S^\star_{\tr_{i+1}}(z,q) - S^\star_{\one+\tr_{i+2}}(z,q) - S^\star_{\one+\tr_{i}}(z,q) + (1-zq)S^\star_{\two+\tr_{i+1}}(z,q) =0, \qquad\! 1 \le i \le \ell-1, \!\!\!\label{eq:finFE_AB2} \\
S^\star_{\zero}(z,q) - 2S^\star_{\one+\ee_\ell}(z,q) + (1-zq) S^\star_{\two}(z,q) = 0. \label{eq:finFE_AB3}
\end{gather}

So now Proposition~\ref{prop:comp1_AB} gives us~\eqref{eq:finFE_AB3} verbatim. Equation~\eqref{eq:finFE_AB1} follows from the fact that $\tr_1=\one+\tr_2$. Additionally, we can use Proposition~\ref{prop:comp2_AB}, together with the facts that $\ee_1+\cdots+\ee_{i}+\tr_{i+1} = \one+\tr_{i+2}$ and $\one+\tr_{i} =\two-\ee_1-\cdots-\ee_{i-1}+\tr_{i+1}$, to obtain~\eqref{eq:finFE_AB2}. Finally, checking the initial conditions ($P^\star_i(0,q)=T^\star_i(0,q)=1$) allows us to complete our proof of Theorem~\ref{thm:mainresult2}.

\section{A bijection with cylindric partitions}\label{sec:6}

Our goal is to demonstrate a bijection between cylindric partitions of profile $(2\ell+1,0)$ on the one hand, with generating function
\[\frac{\bigl(q^1,q^{2\ell+2},q^{2\ell+3};q^{2\ell+3}\bigr)_\infty}{(q;q)_\infty ^2},\]
and pairs of ordinary partitions and $[1,0,\dots,0]$-admissible partitions on $2\ell$ rows on the other, with respective generating functions
\[ \frac{1}{(q;q)_\infty} \qquad \text{and} \qquad \frac{\bigl(q^1,q^{2\ell+2},q^{2\ell+3};q^{2\ell+3}\bigr)_\infty}{(q;q)_\infty}.\] We will follow Proposition 2 of Corteel~\cite{Cort} and its proof. As in Corteel's proof, the first step involves what is essentially ``flattening'' the cylindric partition (meaning to take the two ordinary partitions~$\tau^{(1)}$ and~$\tau^{(2)}$ of the cylindric partition and construct a new ordinary partition~$\tau^{\star}$, where the frequency of all parts~$m$ in~$\tau^{\star}$ is equal to the sum of the frequencies of that part~$m$ in~$\tau^{(1)}$ and~$\tau^{(2)}$) and then computing the conjugate partition of this flattened cylindric partition, and finally assigning colors to the parts of this conjugate. Our color assignment is similar to Corteel's. However, we differ in our method for constructing the pair of partitions $\mu$ and $\nu$ -- in the case $\ell=1$ (corresponding to the Rogers--Ramanujan case, which was the only one considered by Corteel), the map produced by our bijection turns out to be different than that of Corteel.

Throughout this section, when considering colored partitions, we will use subscripts for colors: $i_A$ is a part of size $i$ with color~$A$. Also, we will consider $[1,0,\dots,0]$-admissible partitions on~$2\ell$ rows, and explicitly color their parts, giving parts in rows $2\ell-2A+1$ and $2\ell-2A+2$ the color~$A$. The initial conditions for $[1,0,\dots,0]$-admissible partitions mean that the smallest part that is allowed to appear in row~$j$ is $2\ell-j+2$. Here is an illustration in the case $\ell=3$:
\begin{equation*}
\begin{matrix}
 && && && 7_3 && 9_3 && 11_3 & \\
& && && 6_3 && 8_3 && 10_3 && 12_3\\
 && && 5_2 && 7_2 && 9_2 && 11_2 & \\
& && 4_2 && 6_2 && 8_2 && 10_2 && 12_2\\
 && 3_1 && 5_1 && 7_1 && 9_1 && 11_1 & \\
& 2_1 && 4_1 && 6_1 && 8_1 && 10_1 && 12_1
 \end{matrix}\quad\dots\, .
\end{equation*}
Using this coloring convention, it will be helpful to give an alternative to the downward-path requirement for $[1,0,\dots,0]$-admissible partitions on $2\ell$ rows.

\begin{Proposition} \label{pr:eq}
A colored partition on $2\ell$ rows is a $[1,0,\dots,0]$-admissible partition if and only if it satisfies the following conditions:
\begin{itemize}\itemsep=0pt
\item All parts $i_A$ are distinct.
\item Parts $i_A$ and $i_B$ $(A\ne B)$ cannot both appear.
\item No parts have size $1$.
\item For $1\le i < \ell$, parts of size $2i$ and $2i+1$ can only occur with colors in the set $\{1,\dots,i\}$.
\item Suppose $i< j$. The parts $i_A$ and $j_{B}$ $(A,B>0)$ are forbidden from appearing together if
\begin{itemize}\itemsep=0pt
\item either $i$ is even and $A \le B$, or $i$ is odd and $A < B$, and furthermore $ \bigl\lfloor\frac j2\bigr\rfloor-\left\lfloor\frac i2\right\rfloor\le B-A$, or
\item either $i$ is even and $B< A$, or $i$ is odd and $B\le A$, and furthermore
$ \bigl\lceil\frac j2\bigr\rceil-\bigl\lceil\frac i2\bigr\rceil\le A-B$.
\end{itemize}
\end{itemize}
\end{Proposition}
\begin{proof}
Suppose $\LH$ is a $[1,0,\dots,0]$-admissible partition on $2\ell$ rows, so $k_0=k=1$. Since $k=1$, it follows that all parts are distinct, and $i_A$ and $i_B$ ($A\ne B$) cannot both appear (as otherwise those parts would be in the same downward path). The initial condition $k_0=1$ means that no parts have size 1 and that, for $1\le A < \ell$, parts with color $A$ cannot have size less than $2A$, as otherwise such parts would be in a downward path with the initial condition $\ficone$.
Finally, suppose that parts $i_A$ and $j_B$ (with $i<j$) lie in the same downward path. If the row that $i_A$ lies in is below the row of $j_B$, then, if $i$ is even, then $A \le B$, and, if $i$ is odd, then $A < B$. In this case, we must have $ \bigl\lfloor\frac j2\bigr\rfloor-\bigl\lfloor\frac i2\bigl\rfloor\le B-A$. On the other hand, if the row that $i_A$ lies in is above the row of $j_B$, then, if $i$ is even, then $B<A$, and, if $i$ is odd, then $B \le A$. In this case, we must have $ \bigl\lceil\frac j2\bigr\rceil-\bigl\lceil\frac i2\bigr\rceil\le A-B$. Taking these two cases together verifies the fifth condition. Thus, we have verified all five conditions above.

Now, suppose that $\LH$ satisfies the five conditions above. We will show that it is a $[1,0,\dots,0]$-admissible partition on $2\ell$ rows by showing that no two parts (including the initial condition $k_0=1$) lie on the same downward path. (Note that we can assume that those parts are distinct by the first condition above.)

Suppose that a downward path includes the $k_0=1$ initial condition. In that case, the possible parts in the path are all of the form $i_A$ with $i< 2A$. But those are exactly the ones forbidden by the third and fourth conditions above.

So now, suppose that we have a downward path that does not include the $k_0=1$ initial condition. Suppose that $i_A$ and $j_B$ both lie on this downward path. If $i=j$, then that pair is forbidden by the first condition above, so suppose without loss of generality that $i<j$. As~in the previous case, if $i_A$ is in a lower row than $j_B$ (while sharing the same downward path), then, if $i$ is even, then $A \le B$, and if $i$ is odd, then $A<B$, and furthermore $ \bigl\lfloor\frac j2\bigr\rfloor-\bigl\lfloor\frac i2\bigr\rfloor\le B-A$. Alternatively, if $i_A$ is in a higher row than $j_B$ , then, if $i$ is even, then $B<A$, and if $i$ is odd, then $B \le A$, and furthermore
$ \bigl\lceil\frac j2\bigr\rceil-\bigl\lceil\frac i2\bigr\rceil\le A-B$. But these are exactly the cases forbidden in the fifth condition above. We conclude that there are no downward paths with multiple parts appearing in them, and that $\LH$ is a $[1,0,\dots,0]$-admissible partition on $2\ell$ rows.
\end{proof}

\begin{Theorem} \label{thm:bij}
Fix $\ell\ge1$. For all $n\ge0$, the following are equinumerous:
\begin{itemize}\itemsep=0pt
\item The set $S_1(n)$ of cylindric partitions of weight $n$ with profile $(2\ell+1,0)$.
\item The set $S_2(n)$ of colored partitions of weight $n$ with the following properties:
\begin{itemize}\itemsep=0pt
\item All colors are in the set $\{0,1,\dots,\ell\}$.
\item Any two parts of size $i$ that appear must have the same color.
\item Parts of size $1$ can only occur with color $0$.
\item For $1\le i < \ell$, parts of size $2i$ and $2i+1$ can only occur with colors $\{0,1,\dots,i\}$.
\item Suppose $i<j$. If the parts $i_A$ and $j_{B}$ both appear in $\Lambda$, then both of the following inequalities are true:
\begin{align}
& \biggl\lfloor\frac j2\biggr\rfloor-\biggl\lfloor\frac i2\biggr\rfloor \ge B - A, \label{eq:floor} \\
& \biggl\lceil\frac j2\biggr\rceil-\biggl\lceil\frac i2\biggr\rceil \ge A - B. \label{eq:ceil}
\end{align}
\end{itemize}
\item The set $S_3(n)$ consisting of pairs $(\mu,\nu)$, with the sum of the weights of $\mu$ and $\nu$ is $n$, where $\mu$ is an ordinary partition and $\nu$ is a colored $[1,0,\dots,0]$-admissible partition of $2\ell$ rows $($counted by $P_0(z,q))$.
\end{itemize}
\end{Theorem}
\begin{proof}
We will prove this by exhibiting bijections between $S_1(n)$ and $S_2(n)$ and between $S_2(n)$ and $S_3(n)$, which then gives us a bijection between $S_1(n)$ and $S_3(n)$.

{\it Part $1$: $S_1(n)$ to $S_2(n)$.} First, let us construct the map from $S_1(n)$ to $S_2(n)$. Consider a~cylindric partition $\LH$ of weight $n$ of profile $(2\ell+1,0)$, so $\LH$ is in $S_1(n)$. This corresponds to a~pair of partitions $\bigl(\tau^{(1)},\tau^{(2)}\bigr)$ such that, for all $i$,
\begin{equation*}
\tau_i^{(1)} \ge \tau_i^{(2)} \qquad \text{and}\qquad \tau_i^{(2)} \ge \tau_{i+2\ell+1}^{(1)}.
\end{equation*}
Let $\gamma_1(m)$ and $\gamma_2(m)$ be the number of parts of $\tau^{(1)}$ and $\tau^{(2)}$ (respectively) that are at least $m$. We construct from $\LH$ a colored partition $\Lambda$, where the allowable colors are $0,1,\dots,\ell$. The number of parts of $\Lambda$ will be \smash{$\tau_1^{(1)}$: $\Lambda = \lambda_1+ \cdots + \lambda_{\tau_1^{(1)}}$}. For \smash{$1\le m \le \tau_1^{(1)}$}, $\lambda_m = \gamma_1(m)+\gamma_2(m)$, and $\lambda_m$ has color given by \[c_m=\left\lfloor\frac{\gamma_1(m)-\gamma_2(m)}2\right\rfloor.\] Note that, for all $i$, $0 \le \gamma_1(m)-\gamma_2(m) \le 2\ell+1$.

We claim that $\Lambda$ is in $S_2(n)$. The first condition of $S_2(n)$ is satisfied by construction. Since~$\gamma_1(m)$ and $\gamma_2(m)$ are both nonincreasing with respect to $m$, if $\lambda_{m}=\lambda_{m+1}$, then $\gamma_1(m)=\gamma_1(m+1)$ and $\gamma_2(m)=\gamma_2(m+1)$, and thus the colors are the same:{\samepage
\begin{equation*}
\left\lfloor\frac{\gamma_1(m)-\gamma_2(m)}2\right\rfloor=\left\lfloor\frac{\gamma_1(m+1)-\gamma_2(m+1)}2\right\rfloor.
\end{equation*}
Hence, the second condition of $S_2(n)$ is also satisfied.}

Next, suppose that $\lambda_m=1$ for some $m$. Since $\lambda_1= \gamma_1(m)+\gamma_2(m)$ and $\gamma_1(m)\ge \gamma_2(m)$, it follows that $\gamma_1(m)=1$, $\gamma_2(m)=0$, and
\[ c_m=\left\lfloor\frac{1-0}2\right\rfloor=0,\] which shows that the third condition of $S_2(n)$ is satisfied.

Now, suppose that $\bigl\lfloor\frac{\lambda_m}2\bigr\rfloor=i$ for some $1\le i <\ell$. This implies that $\gamma_1(m)\le 2i+1$, and thus \begin{equation*}
c_m=\left\lfloor\frac{\gamma_1(m)-\gamma_2(m)}2\right\rfloor\le\left\lfloor\frac{2i+1-0}2\right\rfloor=i.
\end{equation*}
This shows that the only possible colors for $\lambda_m$ are $\{0,1,\dots,i\}$, and thus the fourth condition of $S_2(n)$ is satisfied.

Finally, to show that the two inequalities in the fifth condition of $S_2(n)$ are satisfied, assume that $\lambda_{m_1}=i$ and $\lambda_{m_2}=j$, with $i<j$, and so $m_1>m_2$. Note that, for integers $C$ and $D$,
\[ \left\lfloor\frac{C+D}2\right\rfloor-\left\lfloor\frac{C-D}2\right\rfloor=D.\] We can use this fact to show that
\begin{equation*}
\left\lfloor\frac i2\right\rfloor-A=\left\lfloor\frac {\gamma_1(m_1)+\gamma_2(m_1)}2\right\rfloor-\left\lfloor\frac {\gamma_1(m_1)-\gamma_2(m_1)}2\right\rfloor=\gamma_2(m_1),
\end{equation*}
and similarly $\bigl\lfloor\frac j2\bigr\rfloor-B=\gamma_2(m_2)$. Since $m_1>m_2$, $\gamma_2(m_1)\le\gamma_2(m_2)$, and so we conclude
\begin{equation*}
\left\lfloor\frac i2\right\rfloor-A\le\left\lfloor\frac j2\right\rfloor-B.
\end{equation*}
Rearranging this inequality gives~\eqref{eq:floor}. In the same way, for integers $C$ and $D$,
\[ \left\lceil\frac{C+D}2\right\rceil+\left\lfloor\frac{C-D}2\right\rfloor=C.\] Thus,
\begin{equation*}
\left\lceil\frac i2\right\rceil+A=\left\lceil\frac{\gamma_1(m_1)+\gamma_2(m_1)}2\right\rceil+\left\lfloor\frac {\gamma_1(m_1)-\gamma_2(m_1)}2\right\rfloor=\gamma_1(m_1),
\end{equation*}
and similarly $\bigl\lceil\frac j2\bigr\rceil+B=\gamma_1(m_2)$. Now, using $\gamma_1(m_1)\le\gamma_1(m_2)$ gives
\begin{equation*}
\left\lceil\frac j2\right\rceil+B \ge \left\lceil\frac i2\right\rceil+A,
\end{equation*}
which can be rearranged to produce~\eqref{eq:ceil}. Since all five conditions are satisfied, we accordingly conclude that $\Lambda \in S_2(n)$.

{\it Part $2$: $S_2(n)$ to $S_1(n)$.} Now, let us construct the inverse map (from $S_2(n)$ to $S_1(n)$), and let $\Lambda=\lambda_1+\cdots+\lambda_m$ be a colored partition in $S_2(n)$. Let $c_j$ be the color of $\lambda_j$. We will construct a pair of partitions $\bigl(\tau^{(1)},\tau^{(2)}\bigr)$, where $\tau^{(1)}$ will have $\bigl\lceil\frac {\lambda_1}2\bigr\rceil+c_1$ parts and $\tau^{(2)}$ will have $\bigl\lfloor\frac {\lambda_1}2\bigr\rfloor-c_1$ parts. For \smash{$ 1 \le i \le\bigl\lceil\frac {\lambda_1}2\bigr\rceil+c_1$}, set \smash{$\tau^{(1)}_i$} equal to the number of parts $\lambda_j$ such that \smash{$\bigl\lceil\frac {\lambda_j}2\bigr\rceil+c_j\ge i$}, and for \smash{$ 1 \le i \le\bigl\lfloor\frac {\lambda_1}2\bigr\rfloor-c_1$}, set \smash{$\tau^{(2)}_i$} equal to the number of parts $\lambda_j$ such that \smash{$ \bigl\lfloor\frac {\lambda_j}2\bigr\rfloor-c_j\ge i$}.

It is clear by construction that both $\tau^{(1)}$ and $\tau^{(2)}$ are partitions (their parts are in nonincreasing order). Next, for \smash{$ 1 \le i \le\bigl\lfloor\frac {\lambda_1}2\bigr\rfloor-c_1$}, suppose that $\lambda_j$ is a part such that \smash{$\bigl\lfloor\frac {\lambda_j}2\bigr\rfloor-c_j\ge i$} (and thus it is contributing to the count for \smash{$\tau^{(2)}_i$}). Then, $\lambda_j$ also satisfies \smash{$\bigl\lceil\frac {\lambda_j}2\bigr\rceil+c_j\ge i$}, and thus it contributes to the count for \smash{$\tau^{(1)}_i$}. This shows that \smash{$\tau_i^{(1)} \ge \tau_i^{(2)}$}, as needed.

Finally, suppose that $\lambda_j$ is a part such that
\begin{equation*}
\left\lceil\frac {\lambda_j}2\right\rceil+c_j\ge i+2\ell+1
\end{equation*}
(and thus it is contributing to the count for $\tau^{(2)}_{i+2\ell+1}$). Subtracting $2c_j-1$ from both sides of that inequality produces
\begin{equation*}
\left\lceil\frac {\lambda_j}2\right\rceil-1-c_j\ge i+2\ell-2c_j.
\end{equation*}
Recall that the largest possible color is $\ell$, and note that $ \bigl\lceil\frac {\lambda_j}2\bigr\rceil-1\le\bigl\lfloor\frac {\lambda_j}2\bigr\rfloor$. Hence,
\begin{equation*}
\left\lfloor\frac {\lambda_j}2\right\rfloor-1-c_j\ge i,
\end{equation*}
and so $\lambda_j$ contributes to the count for $\tau^{(2)}_i$. This shows that the
final inequality \smash{$\tau_i^{(2)} \ge \tau_{i+2\ell+1}^{(1)}$} is satisfied. We conclude that \smash{$\bigl(\tau^{(1)},\tau^{(2)}\bigr)$} is a cylindric partition of profile $(2\ell+1,0)$.

{\it Part $3$: $S_2(n)$ to $S_3(n)$.} Next, let us construct the map from $S_2(n)$ to $S_3(n)$. Let $\Lambda \in S_2(n)$. We will construct a pair of partitions $\mu$ and $\nu$, where $\mu$ will be an ordinary partition, while $\nu$ will be a colored partition with colors in the set $\{1,\dots,\ell\}$. (When we add parts to $\mu$, we ``forget'' their colors.) First, we take each part of $\Lambda$ that occurs $t>1$ times, and $t-1$ occurrences of that part are added to $\mu$. Additionally, all parts colored 0 are added to $\mu$. Now carry out the following rules:
\begin{itemize}\itemsep=0pt
\item[$(R_1)$] If $i_{A}$ and $j_{B}$ (with $i<j$) both appear as parts in $\Lambda$, where either $i$ is even and $0<A \le B$, or $i$ is odd and $0 < A < B$, and furthermore
\smash{$\bigl\lfloor\frac j2\bigr\rfloor-\bigl\lfloor\frac i2\bigr\rfloor=B-A$}, then insert $i$ into $\mu$.
\item[$(R_2$)] If $i_{A}$ and $j_{B}$ (with $i<j$) both appear as parts in $\Lambda$, where either $i$ is even and $0<B< A$, or $i$ is odd and $0<B\le A$, and furthermore
$\bigl\lceil\frac j2\bigr\rceil-\bigl\lceil\frac i2\bigr\rceil=A-B$, then insert $j$ into $\mu$.
\end{itemize}
At the end, we form $\nu$ by taking all of the leftover parts from $\Lambda$ not inserted into $\mu$, keeping their colors intact.

Now, $\mu$ is clearly an (ordinary) partition. We claim that $\nu$ is a colored $[1,0,\dots,0]$-admissible partition on $2\ell$ rows through checking the conditions of Proposition~\ref{pr:eq}. The first condition above is met by the fact that duplicated parts were sent to $\mu$, while the second, third, and fourth are met through the properties of $\Lambda$. (Recall that all parts of $\Lambda$ colored 0 are sent to $\nu$.) Finally, the fifth condition nearly follows from the properties of $\Lambda$ expressed in~\eqref{eq:floor} and~\eqref{eq:ceil}, except for the extra cases
\begin{itemize}\itemsep=0pt
\item either $i$ is even and $A \le B$, or $i$ is odd and $A < B$, and furthermore $ \bigl\lfloor\frac j2\bigr\rfloor-\bigl\lfloor\frac i2\bigr\rfloor= B-A$, or
\item either $i$ is even and $B< A$, or $i$ is odd and $B\le A$, and furthermore
$\bigl\lceil\frac j2\bigr\rceil-\bigl\lceil\frac i2\bigr\rceil= A-B$.
\end{itemize}
But, these cases are exactly the ones that were dealt with in the rules $R_1$ and $R_2$, and so we conclude that $\nu$ is, in fact, a $[1,0,\dots,0]$-admissible partition on $2\ell$ rows (with bivariate generating function $P_0(z,q)$).

{\it Part $4$: $S_3(n)$ to $S_2(n)$.} Finally, let us construct the inverse map (from $S_3(n)$ to $S_2(n)$). Let $(\mu,\nu) \in S_3(n)$. We transform this pair into a colored partition $\Lambda$ by carrying out the following operations:
\begin{itemize}\itemsep=0pt
\item[$(R_1')$] If $i$ is in $\mu$, $j_B$ is in $\nu$, $i<j$, and either $i$ is even and $j-i\le2B-1$, or $i$ is odd and $j-i<2B-1$, then insert $i_C$ into $\Lambda$, where $C=B-\bigl\lfloor\frac j2\bigr\rfloor+\bigl\lfloor\frac i2\bigr\rfloor$.
\item[$(R_2')$] If $j$ is in $\mu$, $i_A$ is in $\nu$, $i<j$, and either $i$ is even and $j-i<2A-1$, or $i$ is odd and $j-i\le 2A-1$, then insert $j_C$ into $\Lambda$, where $ C=A+\bigl\lceil\frac j2\bigr\rceil-\bigl\lceil\frac i2\bigr\rceil$.
\end{itemize}
All parts of $\nu$ are then inserted into $\Lambda$ with their colors intact, and finally all remaining parts~$i$ of~$\mu$ are inserted into $\Lambda$. If there already is a part $i_A$ in $\Lambda$, then those parts $i$ are also given color~$A$; otherwise, we insert $i_0$ into $\Lambda$.

We claim that this is the inverse of the map given in Part~3. First, note that neither the map in Part 3 nor this map changes the size of any part -- only colors are (potentially) changed, deleted, or added. Parts that are in $\nu$ do not have their colors changed by either the Part~3 map or this map. If some part of $\Lambda$ that had nonzero color occurred $t$ times and had $t-1$ copies sent to $\mu$ (thus losing their color), those parts have their color restored by the ``if there already is a~part $i_A$ in $\Lambda$, then those parts~$i$ are also given color~$A$'' part of the map. Finally, parts of~$\Lambda$ that had nonzero color and were sent to $\mu$ by either $R_1$ or $R_2$ in the Part~3 map (thus losing their color) are exactly the ones that have their colors restored by $R_1'$ and $R_2'$ in this present map.

For example, suppose that $i$ is even and $i_A$ and $j_B$ are parts of $\Lambda$ that satisfy the conditions of $R_1$ in Part 3, so $i$ is inserted into $\mu$. Since $\bigl\lfloor\frac j2\bigr\rfloor-\bigl\lfloor\frac i2\bigr\rfloor=B-A$, it follows that
\begin{equation*}
2A-1+2\left\lfloor\frac j2\right\rfloor-2\left\lfloor\frac i2\right\rfloor=2B-1.
\end{equation*}
We assumed $i$ is even, and so $ 2\bigl\lfloor\frac i2\bigr\rfloor=i$, and since $A>0$, $2A-1+2\bigl\lfloor\frac j2\bigr\rfloor\ge j-1$, so $j-i\le 2B-1$. Hence, the hypotheses of $R_1'$ are met, and we insert $i_C$ into $\Lambda$, where $ C=B-\bigl\lfloor\frac j2\bigr\rfloor+\bigl\lfloor\frac i2\bigr\rfloor$, but then it is easy to see that this color $C$ matches exactly the original color $A$, and so the correct color has been restored to $i$ in $\Lambda$. The other cases ($R_1$ and $i$ odd, along with $R_2$ and both parities of $i$) follow similarly.

We thus conclude that we have successfully produced bijections between $S_1(n)$ and $S_2(n)$ and between $S_2(n)$ and $S_3(n)$. This gives us, as desired, a bijection between cylindric partitions of profile $(2\ell+1,0)$ and pairs $(\mu,\nu)$, where $\mu$ is an ordinary partition and $\nu$ is a $[1,0,\dots,0]$-admissible partition on $2\ell$ rows.
\end{proof}

Let us present an example to demonstrate our work. Suppose that $\ell=3$, and our original cylindric partition \smash{$\LH=\bigl(\tau^{(1)},\tau^{(2)}\bigr)$} has
$\tau^{(1)}=9+8+8+5+5+3+3+3+3+3$,
$\tau^{(2)}=9+4+4+4+2$, so $\LH \in S_1(73)$.
Then,
\smash{$\Lambda=15_2+15_2+14_3+9_0+6_2+4_1+4_1+4_1+2_0$}.
Observe that $\Lambda$ meets the properties to be in $S_2(73)$.

We now construct $\mu$ and $\nu$. Since $15_2$ and $4_1$ occur multiple times in $\Lambda$, we place all but one occurrence of each into $\mu$ (with colors deleted), and we also place the parts colored 0 (9 and 2) into $\mu,$ so right now, we have $\mu = 15+9+4+4+2.$ Following rule $R_2$ above for $15_2+14_3$, we see that we insert the final copy of 15 into $\mu$, and following rule $R_1$ for $6_2+4_1$, we see that we also insert the final copy of 4 into $\mu$. This leaves us with $\mu = 15+15+9+4+4+4+2$ and $\nu = 14_3+6_2$, where $\mu$ is an ordinary partition and $\nu$ is a $[1,0,0,0]$-admissible partition.

The above work was for $[1,0,\dots,0]$-admissible partitions on $2\ell$ rows. For the case of admissible partitions on $2\ell-1$ rows, we now have a bijection between cylindric partitions of profile $(2\ell,0)$ and pairs of ordinary partitions and $[1,0,\dots,0]$-admissible partitions on $2\ell-1$ rows. It is straightforward to adapt the proof of Theorem~\ref{thm:bij}. The main difference is that, in the colored partition $\Lambda \in S_2(n)$ that is constructed, parts with color $\ell$ cannot have odd size, which parallels exactly how parts with color $\ell$ in the admissible partition $\nu$ that is constructed must have even size (again shown for $\ell=3$):
\begin{equation*}
\begin{matrix}
& && && 6_3 && 8_3 && 10_3 && 12_3\\
 && && 5_2 && 7_2 && 9_2 && 11_2 & \\
& && 4_2 && 6_2 && 8_2 && 10_2 && 12_2\\
 && 3_1 && 5_1 && 7_1 && 9_1 && 11_1 & \\
& 2_1 && 4_1 && 6_1 && 8_1 && 10_1 && 12_1
 \end{matrix}\quad\dots\,.
\end{equation*}

\section{Conclusion}\label{sec:conc7}

We conclude by listing some questions suggested by the present research. We focus on the case with an even number of rows (Theorem~\ref{thm:mainresult1}) corresponding to the Andrews--Gordon identities, but many of the comments below would equally apply to the other case.
\begin{itemize}\itemsep=0pt
 \item As we have seen, taking the usual bivariate Andrews--Gordon and Andrews--Bressoud multisums, which have elegant combinatorial interpretations, and altering the $x^{N_1+\cdots+N_k}$ factor to $z^{N_1}$ produces new bivariate multisums that still have elegant combinatorial interpretations. Can the same be done for other families of multisum identities? One natural candidate would be the G\"ollnitz--Gordon--Andrews identities~\cite{And_GGA1,And_GGA2}.
 \end{itemize}

\begin{Theorem}
Fix $\ell\ge 1$, and let $0\le i\le \ell$. Let $G_{\ell,i,m}(n)$ denote the number of partitions of $n$ with exactly $m$ parts satisfying, for all $j\ge1$, $f_{2j-1}\le 1$ and
$f_{2j}+f_{2j+1}+f_{2j+2}\le\ell$, and also satisfying $f_1+f_2\le i$. Then,
\begin{gather*}
\sum_{m,n \ge 0} G_{\ell,i,m}(n) x^m q^n = \sum_{n_1,n_2,\dots,n_{\ell}\ge 0} \frac{x^{N_1+\cdots+N_\ell} \bigl(-q,q^2\bigr)_{N_1}q^{N_1^2+2(N_2^2+\cdots+N_{\ell}^2+N_{i+1}+N_{i+2}+\cdots+N_{\ell})}}{\bigl(q^2;q^2\bigr)_{n_1}\bigl(q^2;q^2\bigr)_{n_2}\cdots\bigl(q^2;q^2\bigr)_{n_{\ell}}}
\end{gather*}
and, if $x=1$ in the above equation, we have the infinite product
\begin{align*}
&\sum_{n_1,n_2,\dots,n_{\ell}\ge 0} \frac{ \bigl(-q,q^2\bigr)_{N_1}q^{N_1^2+2(N_2^2+\cdots+N_{\ell}^2+N_{i+1}+N_{i+2}+\cdots+N_{\ell})}}{\bigl(q^2;q^2\bigr)_{n_1}\bigl(q^2;q^2\bigr)_{n_2}\cdots\bigl(q^2;q^2\bigr)_{n_{\ell}}} \\
&\qquad= \frac{\bigl(q^2;q^4\bigr)_\infty\bigl(q^{2i+1},q^{4\ell-2i+3},q^{4\ell+4};q^{4\ell+4}\bigr)_\infty}{(q;q)_\infty}.
\end{align*}
\end{Theorem}
So, are there some combinatorial objects (that feel similar to the colored partitions of this paper) with bivariate generating function
\begin{equation*}
\sum_{n_1,n_2,\dots,n_{\ell}\ge 0} \frac{z^{N_1} \bigl(-q,q^2\bigr)_{N_1}q^{N_1^2+2(N_2^2+\cdots+N_{\ell}^2+N_{i+1}+N_{i+2}+\cdots+N_{\ell})}}{\bigl(q^2;q^2\bigr)_{n_1}\bigl(q^2;q^2\bigr)_{n_2}\cdots\bigl(q^2;q^2\bigr)_{n_{\ell}}}?
\end{equation*}
\begin{itemize}\itemsep=0pt
\item Can the ideas and results of this paper be extended to obtain a multisum for the $k>1,$ $\ell>1$ cases of the conjectures of Capparelli, Meurman, Primc, and Primc? To be more specific, let us restrict ourselves to the case with $2\ell$ rows. We have a bivariate generating function (Theorem~\ref{thm:mainresult1}) in the case $k=1$ and $\ell \in \NN$. But, as observed in~\cite{CMPP}, the case $\ell=1$ and $k\in\NN$ also ends up literally as Theorem~\ref{thm:AndGor}, with bivariate generating function~\eqref{eq:AndGorxq}. Can these ideas be combined appropriately to give a multisum for general $k$, $\ell$? Right now, even the case $k=2$ and $\ell=2$ is open.

\item There exist many refinements and modifications of the Andrews--Gordon identities. For example, the parity of parts has been investigated~\cite{AndParity,KimYee,KurThesis,KurParity}. A typical case to study here is partitions counted by $B_{\ell,i}(n)$ in Theorem~\ref{thm:Gordon} that further have the requirement that even parts must appear an even number of times. Once again, there is a bivariate multisum for these partitions with a factor of $x^{N_1+\cdots+N_k}$; if this is modified to $z^{N_1}$, is there a nice subfamily of these generalized partitions that are now counted?

\item Another combinatorial interpretation of the sum on the right side of~\eqref{eq:AndGorxq} involves Durfee squares and rectangles~\cite{AndDurf} (see also~\cite{ADJM}). One nice feature of this combinatorial interpretation is that the refinement obtained provides information about the values of each of the variables $n_1,\dots,n_\ell$ in the sum. Again, is there a corresponding combinatorial interpretation for the bivariate generating functions $P_i(z,q)$ and $P_i^\star(z,q)$ of the Capparelli--Meurman--Primc--Primc colored partitions?
\end{itemize}

\appendix

\section{Appendix: Explicit computations}\label{sec:appendix}
In this appendix, we include explicit examples of some of our computations in the case $\ell=4$ in order to illustrate our proofs.

\subsection[Theorem~\ref{thm:funceq1A} for l=4]{Theorem~\ref{thm:funceq1A} for $\boldsymbol{\ell=4}$}\label{subsec:combspec}

Here, we have the following initial picture, where exactly one of the $k_i$ equals 1, and the rest are zero:
\begin{equation*}
\begin{matrix}
k_1 && f_1 && f_3 && f_5 && f_7 && f_9 && f_{11} && f_{13} & \\
& \cdot && f_2 && f_4 && f_6 && f_8 && f_{10} && f_{12} && f_{14} \\
k_3 && f_1 && f_3 && f_5 && f_7 && f_9 && f_{11} && f_{13} & \\
& \cdot && f_2 && f_4 && f_6 && f_8 && f_{10} && f_{12} && f_{14} \\
k_4 && f_1 && f_3 && f_5 && f_7 && f_9 && f_{11} && f_{13} & \\
& \cdot && f_2 && f_4 && f_6 && f_8 && f_{10} && f_{12} && f_{14} \\
k_2 && f_1 && f_3 && f_5 && f_7 && f_9 && f_{11} && f_{13} & \\
& k_0 && f_2 && f_4 && f_6 && f_8 && f_{10} && f_{12} && f_{14}
 \end{matrix}\quad\dots\,.
\end{equation*}

We display frequency arrays for the five generating functions that we care about, with their initial conditions shown as $\ficone$, and $\cdot$ representing parts that are forbidden from appearing:\looseness=1
\begin{gather}
P_4(z,q)\colon \quad
\begin{matrix}
\cdot && \cdot && \cdot && f_5 && f_7 && f_9 && f_{11} && f_{13} & \\
& \cdot && \cdot && f_4 && f_6 && f_8 && f_{10} && f_{12} && f_{14} \\
\cdot && \cdot && f_3 && f_5 && f_7 && f_9 && f_{11} && f_{13} & \\
& \cdot && f_2 && f_4 && f_6 && f_8 && f_{10} && f_{12} && f_{14} \\
\ficone && f_1 && f_3 && f_5 && f_7 && f_9 && f_{11} && f_{13} & \\
& \cdot && f_2 && f_4 && f_6 && f_8 && f_{10} && f_{12} && f_{14} \\
\cdot && \cdot && f_3 && f_5 && f_7 && f_9 && f_{11} && f_{13} & \\
& \cdot && \cdot && f_4 && f_6 && f_8 && f_{10} && f_{12} && f_{14} \vspace{3mm}
 \end{matrix}\quad\dots\,,
\nonumber\\
P_3(z,q)\colon \quad
\begin{matrix}
\cdot && \cdot && f_3 && f_5 && f_7 && f_9 && f_{11} && f_{13} & \\
& \cdot && f_2 && f_4 && f_6 && f_8 && f_{10} && f_{12} && f_{14} \\
\ficone && f_1 && f_3 && f_5 && f_7 && f_9 && f_{11} && f_{13} & \\
& \cdot && f_2 && f_4 && f_6 && f_8 && f_{10} && f_{12} && f_{14} \\
\cdot && \cdot && f_3 && f_5 && f_7 && f_9 && f_{11} && f_{13} & \\
& \cdot && \cdot && f_4 && f_6 && f_8 && f_{10} && f_{12} && f_{14} \\
\cdot && \cdot && \cdot && f_5 && f_7 && f_9 && f_{11} && f_{13} & \\
& \cdot && \cdot && \cdot && f_6 && f_8 && f_{10} && f_{12} && f_{14} \vspace{3mm}
 \end{matrix}\quad\dots\,,
\nonumber\\
P_2(z,q)\colon \quad
\begin{matrix}
\cdot && \cdot && \cdot && \cdot && f_7 && f_9 && f_{11} && f_{13} & \\
& \cdot && \cdot && \cdot && f_6 && f_8 && f_{10} && f_{12} && f_{14} \\
\cdot && \cdot && \cdot && f_5 && f_7 && f_9 && f_{11} && f_{13} & \\
& \cdot && \cdot && f_4 && f_6 && f_8 && f_{10} && f_{12} && f_{14} \\
\cdot && \cdot && f_3 && f_5 && f_7 && f_9 && f_{11} && f_{13} & \\
& \cdot && f_2 && f_4 && f_6 && f_8 && f_{10} && f_{12} && f_{14} \\
\ficone && f_1 && f_3 && f_5 && f_7 && f_9 && f_{11} && f_{13} & \\
& \cdot && f_2 && f_4 && f_6 && f_8 && f_{10} && f_{12} && f_{14} \vspace{3mm}
 \end{matrix}\quad\dots\,,
\nonumber\\
P_1(z,q)\colon \quad
\begin{matrix}
\ficone && f_1 && f_3 && f_5 && f_7 && f_9 && f_{11} && f_{13} & \\
& \cdot && f_2 && f_4 && f_6 && f_8 && f_{10} && f_{12} && f_{14} \\
\cdot && \cdot && f_3 && f_5 && f_7 && f_9 && f_{11} && f_{13} & \\
& \cdot && \cdot && f_4 && f_6 && f_8 && f_{10} && f_{12} && f_{14} \\
\cdot && \cdot && \cdot && f_5 && f_7 && f_9 && f_{11} && f_{13} & \\
& \cdot && \cdot && \cdot && f_6 && f_8 && f_{10} && f_{12} && f_{14} \\
\cdot && \cdot && \cdot && \cdot && f_7 && f_9 && f_{11} && f_{13} & \\
& \cdot && \cdot && \cdot && \cdot && f_8 && f_{10} && f_{12} && f_{14}\vspace{3mm}
 \end{matrix}\quad\dots\,,
\nonumber\\
P_0(z,q)\colon \quad
\begin{matrix}
\cdot && \cdot && \cdot && \cdot && \cdot && f_9 && f_{11} && f_{13} & \\
& \cdot && \cdot && \cdot && \cdot && f_8 && f_{10} && f_{12} && f_{14} \\
\cdot && \cdot && \cdot && \cdot && f_7 && f_9 && f_{11} && f_{13} & \\
& \cdot && \cdot && \cdot && f_6 && f_8 && f_{10} && f_{12} && f_{14} \\
\cdot && \cdot && \cdot && f_5 && f_7 && f_9 && f_{11} && f_{13} & \\
& \cdot && \cdot && f_4 && f_6 && f_8 && f_{10} && f_{12} && f_{14} \\
\cdot && \cdot && f_3 && f_5 && f_7 && f_9 && f_{11} && f_{13} & \\
& \ficone && f_2 && f_4 && f_6 && f_8 && f_{10} && f_{12} && f_{14} \\
 \end{matrix}\quad\dots\, .\!\!\!\label{eq:a3}
\end{gather}

We now show how to deduce the functional equations corresponding to~\eqref{eq:origFE_gf1} and~\eqref{eq:origFE_gf2} in this particular case:
\begin{itemize}\itemsep=0pt
 \item $P_0(z,q) = P_1(zq,q)$: This is clear by inspection: simply ``flip'' the picture in~\eqref{eq:a3} upside down. This flipping trick will be used constantly throughout the rest of this appendix.\looseness=1
\item $P_1(z,q) =
zq^1 P_{1}\bigl(zq^{2},q\bigr) + P_{2}(zq,q)$: Consider a partition counted by $P_1(z,q)$. Either there is a 1 in this partition, or there is not. If there is, we have the following picture:\looseness=1
\begin{equation*}
\begin{matrix}
\cdot && {\color{red}{1}} && f_3 && f_5 && f_7 && f_9 && f_{11} && f_{13} & \\
& \cdot && \rs && f_4 && f_6 && f_8 && f_{10} && f_{12} && f_{14} \\
\cdot && \cdot && \rs && f_5 && f_7 && f_9 && f_{11} && f_{13} & \\
& \cdot && \cdot && \rs && f_6 && f_8 && f_{10} && f_{12} && f_{14} \\
\cdot && \cdot && \cdot && \rs && f_7 && f_9 && f_{11} && f_{13} & \\
& \cdot && \cdot && \cdot && \rs && f_8 && f_{10} && f_{12} && f_{14} \\
\cdot && \cdot && \cdot && \cdot && \rs && f_9 && f_{11} && f_{13} & \\
& \cdot && \cdot && \cdot && \cdot && \rs && f_{10} && f_{12} && f_{14}
 \end{matrix}\quad\dots\,.
\end{equation*}
Here (and following), the red $\color{red}{1}$ indicates the location of the 1 we are assuming to be in this partition, while the red $\rs$s show the frequencies that are forced to equal zero by the inclusion of this part. We can see that this corresponds to $zqP_1\bigl(zq^2,q\bigr)$. The other possibility is that there is no 1
\begin{equation*}
\begin{matrix}
\cdot && \rs && f_3 && f_5 && f_7 && f_9 && f_{11} && f_{13} & \\
& \cdot && f_2 && f_4 && f_6 && f_8 && f_{10} && f_{12} && f_{14} \\
\cdot && \cdot && f_3 && f_5 && f_7 && f_9 && f_{11} && f_{13} & \\
& \cdot && \cdot && f_4 && f_6 && f_8 && f_{10} && f_{12} && f_{14} \\
\cdot && \cdot && \cdot && f_5 && f_7 && f_9 && f_{11} && f_{13} & \\
& \cdot && \cdot && \cdot && f_6 && f_8 && f_{10} && f_{12} && f_{14} \\
\cdot && \cdot && \cdot && \cdot && f_7 && f_9 && f_{11} && f_{13} & \\
& \cdot && \cdot && \cdot && \cdot && f_8 && f_{10} && f_{12} && f_{14}
 \end{matrix}\quad\dots\,,
\end{equation*}
which corresponds to $P_2(zq,q)$ by flipping the array.

\item $P_2(z,q) = zq^1 P_{2}\bigl(zq^{2},q\bigr) +
zq^2 P_{1}\bigl(zq^{3},q\bigr) + P_3(zq,q)$: Either there is a 1
\begin{equation*}
\begin{matrix}
\cdot && \cdot && \cdot && \cdot && \rs && f_9 && f_{11} && f_{13} & \\
& \cdot && \cdot && \cdot && \rs && f_8 && f_{10} && f_{12} && f_{14} \\
\cdot && \cdot && \cdot && \rs && f_7 && f_9 && f_{11} && f_{13} & \\
& \cdot && \cdot && \rs && f_6 && f_8 && f_{10} && f_{12} && f_{14} \\
\cdot && \cdot && \rs && f_5 && f_7 && f_9 && f_{11} && f_{13} & \\
& \cdot && \rs && f_4 && f_6 && f_8 && f_{10} && f_{12} && f_{14} \\
\cdot && \color{red}{1} && f_3 && f_5 && f_7 && f_9 && f_{11} && f_{13} & \\
& \cdot && \rs && f_4 && f_6 && f_8 && f_{10} && f_{12} && f_{14}
 \end{matrix}\quad\dots\,,
\end{equation*}
which corresponds to $zq^1 P_{2}\bigl(zq^{2},q\bigr)$, or there is a 2 in the final row\footnote{In a small abuse of notation, we are using {\color{red}{2}} to indicate $f_2=1$ here (and similarly throughout the rest of the appendix).}
\begin{equation*}
\begin{matrix}
\cdot && \cdot && \cdot && \cdot && \rs && \rs && f_{11} && f_{13} & \\
& \cdot && \cdot && \cdot && \rs && \rs && f_{10} && f_{12} && f_{14} \\
\cdot && \cdot && \cdot && \rs && \rs && f_9 && f_{11} && f_{13} & \\
& \cdot && \cdot && \rs && \rs && f_8 && f_{10} && f_{12} && f_{14} \\
\cdot && \cdot && \rs && \rs && f_7 && f_9 && f_{11} && f_{13} & \\
& \cdot && \rs && \rs && f_6 && f_8 && f_{10} && f_{12} && f_{14} \\
\cdot && \rs && \rs && f_5 && f_7 && f_9 && f_{11} && f_{13} & \\
& \cdot && \color{red}{2} && f_4 && f_6 && f_8 && f_{10} && f_{12} && f_{14}
 \end{matrix}\quad\dots\,,
\end{equation*}
which corresponds to $zq^2 P_{1}\bigl(zq^{3},q\bigr)$, or neither of these occur
\begin{equation*}
\begin{matrix}
\cdot && \cdot && \cdot && \cdot && f_7 && f_9 && f_{11} && f_{13} & \\
& \cdot && \cdot && \cdot && f_6 && f_8 && f_{10} && f_{12} && f_{14} \\
\cdot && \cdot && \cdot && f_5 && f_7 && f_9 && f_{11} && f_{13} & \\
& \cdot && \cdot && f_4 && f_6 && f_8 && f_{10} && f_{12} && f_{14} \\
\cdot && \cdot && f_3 && f_5 && f_7 && f_9 && f_{11} && f_{13} & \\
& \cdot && f_2 && f_4 && f_6 && f_8 && f_{10} && f_{12} && f_{14} \\
\cdot && \rs && f_3 && f_5 && f_7 && f_9 && f_{11} && f_{13} & \\
& \cdot && \rs && f_4 && f_6 && f_8 && f_{10} && f_{12} && f_{14}
 \end{matrix}\quad\dots\,,
\end{equation*}
which corresponds to $P_3(zq,q)$.
\item $P_3(z,q) =
zq^1 P_{3}\bigl(zq^{2},q\bigr) +
zq^2 P_{2}\bigl(zq^{3},q\bigr) +
zq^3 P_{1}\bigl(zq^{4},q\bigr) + P_{4}(zq,q)$: Either there is a~1
\begin{equation*}
\begin{matrix}
\cdot && \cdot && \rs && f_5 && f_7 && f_9 && f_{11} && f_{13} & \\
& \cdot && \rs && f_4 && f_6 && f_8 && f_{10} && f_{12} && f_{14} \\
\cdot && \color{red}{1} && f_3 && f_5 && f_7 && f_9 && f_{11} && f_{13} & \\
& \cdot && \rs && f_4 && f_6 && f_8 && f_{10} && f_{12} && f_{14} \\
\cdot && \cdot && \rs && f_5 && f_7 && f_9 && f_{11} && f_{13} & \\
& \cdot && \cdot && \rs && f_6 && f_8 && f_{10} && f_{12} && f_{14} \\
\cdot && \cdot && \cdot && \rs && f_7 && f_9 && f_{11} && f_{13} & \\
& \cdot && \cdot && \cdot && \rs && f_8 && f_{10} && f_{12} && f_{14}
 \end{matrix}\quad\dots\,,
\end{equation*}
which corresponds to $zq^1 P_{3}\bigl(zq^{2},q\bigr)$, or there is a 2 in the second row
\begin{equation*}
\begin{matrix}
\cdot && \cdot && \rs && f_5 && f_7 && f_9 && f_{11} && f_{13} & \\
& \cdot && \color{red}{2} && f_4 && f_6 && f_8 && f_{10} && f_{12} && f_{14} \\
\cdot && \rs && \rs && f_5 && f_7 && f_9 && f_{11} && f_{13} & \\
& \cdot && \rs && \rs && f_6 && f_8 && f_{10} && f_{12} && f_{14} \\
\cdot && \cdot && \rs && \rs && f_7 && f_9 && f_{11} && f_{13} & \\
& \cdot && \cdot && \rs && \rs && f_8 && f_{10} && f_{12} && f_{14} \\
\cdot && \cdot && \cdot && \rs && \rs && f_9 && f_{11} && f_{13} & \\
& \cdot && \cdot && \cdot && \rs && \rs && f_{10} && f_{12} && f_{14}
 \end{matrix}\quad\dots\,,
\end{equation*}
which corresponds to $zq^2 P_{2}\bigl(zq^{3},q\bigr)$, or there is a 3 in the first row
\begin{equation*}
\begin{matrix}
\cdot && \cdot && \color{red}{3} && f_5 && f_7 && f_9 && f_{11} && f_{13} & \\
& \cdot && \rs && \rs && f_6 && f_8 && f_{10} && f_{12} && f_{14} \\
\cdot && \rs && \rs && \rs && f_7 && f_9 && f_{11} && f_{13} & \\
& \cdot && \rs && \rs && \rs && f_8 && f_{10} && f_{12} && f_{14} \\
\cdot && \cdot && \rs && \rs && \rs && f_9 && f_{11} && f_{13} & \\
& \cdot && \cdot && \rs && \rs && \rs && f_{10} && f_{12} && f_{14} \\
\cdot && \cdot && \cdot && \rs && \rs && \rs && f_{11} && f_{13} & \\
& \cdot && \cdot && \cdot && \rs && \rs && \rs && f_{12} && f_{14}
 \end{matrix}\quad\dots\,,
\end{equation*}
which corresponds to $zq^3 P_{1}\bigl(zq^{4},q\bigr)$, or none of these occur
\begin{equation*}
\begin{matrix}
\cdot && \cdot && \rs && f_5 && f_7 && f_9 && f_{11} && f_{13} & \\
& \cdot && \rs && f_4 && f_6 && f_8 && f_{10} && f_{12} && f_{14} \\
\cdot && \rs && f_3 && f_5 && f_7 && f_9 && f_{11} && f_{13} & \\
& \cdot && f_2 && f_4 && f_6 && f_8 && f_{10} && f_{12} && f_{14} \\
\cdot && \cdot && f_3 && f_5 && f_7 && f_9 && f_{11} && f_{13} & \\
& \cdot && \cdot && f_4 && f_6 && f_8 && f_{10} && f_{12} && f_{14} \\
\cdot && \cdot && \cdot && f_5 && f_7 && f_9 && f_{11} && f_{13} & \\
& \cdot && \cdot && \cdot && f_6 && f_8 && f_{10} && f_{12} && f_{14}
 \end{matrix}\quad\dots\,,
\end{equation*}
which corresponds to $P_{4}(zq,q)$.

\item $P_4(z,q) =
zq^1 P_{4}\bigl(zq^{2},q\bigr) +
zq^2 P_{3}\bigl(zq^{3},q\bigr) +
zq^3 P_{2}\bigl(zq^{4},q\bigr) +
zq^4 P_{1}\bigl(zq^{5},q\bigr) + P_4(zq,q)$.

Either a 1 occurs in the partition
\begin{equation*}
\begin{matrix}
\cdot && \cdot && \cdot && \rs && f_7 && f_9 && f_{11} && f_{13} & \\
& \cdot && \cdot && \rs && f_6 && f_8 && f_{10} && f_{12} && f_{14} \\
\cdot && \cdot && \rs && f_5 && f_7 && f_9 && f_{11} && f_{13} & \\
& \cdot && \rs && f_4 && f_6 && f_8 && f_{10} && f_{12} && f_{14} \\
\cdot && \color{red}{1} && f_3 && f_5 && f_7 && f_9 && f_{11} && f_{13} & \\
& \cdot && \rs && f_4 && f_6 && f_8 && f_{10} && f_{12} && f_{14} \\
\cdot && \cdot && \rs && f_5 && f_7 && f_9 && f_{11} && f_{13} & \\
& \cdot && \cdot && \rs && f_6 && f_8 && f_{10} && f_{12} && f_{14}
 \end{matrix}\quad\dots\,,
\end{equation*}
corresponding to $zq^1 P_{4}\bigl(zq^{2},q\bigr)$, or a 2 occurs in the sixth row
\begin{equation*}
\begin{matrix}
\cdot && \cdot && \cdot && \rs && \rs && f_9 && f_{11} && f_{13} & \\
& \cdot && \cdot && \rs && \rs && f_8 && f_{10} && f_{12} && f_{14} \\
\cdot && \cdot && \rs && \rs && f_7 && f_9 && f_{11} && f_{13} & \\
& \cdot && \rs && \rs && f_6 && f_8 && f_{10} && f_{12} && f_{14} \\
\cdot && \rs && \rs && f_5 && f_7 && f_9 && f_{11} && f_{13} & \\
& \cdot && \color{red}{2} && f_4 && f_6 && f_8 && f_{10} && f_{12} && f_{14} \\
\cdot && \cdot && \rs && f_5 && f_7 && f_9 && f_{11} && f_{13} & \\
& \cdot && \cdot && \rs && f_6 && f_8 && f_{10} && f_{12} && f_{14}
 \end{matrix}\quad\dots\,,
\end{equation*}
corresponding to $zq^2 P_{3}\bigl(zq^{3},q\bigr),$ or a 3 occurs in the seventh row
\begin{equation*}
\begin{matrix}
\cdot && \cdot && \cdot && \rs && \rs && \rs && f_{11} && f_{13} & \\
& \cdot && \cdot && \rs && \rs && \rs && f_{10} && f_{12} && f_{14} \\
\cdot && \cdot && \rs && \rs && \rs && f_9 && f_{11} && f_{13} & \\
& \cdot && \rs && \rs && \rs && f_8 && f_{10} && f_{12} && f_{14} \\
\cdot && \rs && \rs && \rs && f_7 && f_9 && f_{11} && f_{13} & \\
& \cdot && \rs && \rs && f_6 && f_8 && f_{10} && f_{12} && f_{14} \\
\cdot && \cdot && \color{red}{3} && f_5 && f_7 && f_9 && f_{11} && f_{13} & \\
& \cdot && \cdot && \rs && f_6 && f_8 && f_{10} && f_{12} && f_{14}
 \end{matrix}\quad\dots\,,
\end{equation*}
corresponding to $zq^3 P_{2}\bigl(zq^{4},q\bigr)$, or a 4 occurs in the eighth row
\begin{equation*}
\begin{matrix}
\cdot && \cdot && \cdot && \rs && \rs && \rs && \rs && f_{13} & \\
& \cdot && \cdot && \rs && \rs && \rs && \rs && f_{12} && f_{14} \\
\cdot && \cdot && \rs && \rs && \rs && \rs && f_{11} && f_{13} & \\
& \cdot && \rs && \rs && \rs && \rs && f_{10} && f_{12} && f_{14} \\
\cdot && \rs && \rs && \rs && \rs && f_9 && f_{11} && f_{13} & \\
& \cdot && \rs && \rs && \rs && f_8 && f_{10} && f_{12} && f_{14} \\
\cdot && \cdot && \rs && \rs && f_7 && f_9 && f_{11} && f_{13} & \\
& \cdot && \cdot && \color{red}{4} && f_6 && f_8 && f_{10} && f_{12} && f_{14} \\
 \end{matrix}\quad\dots,
\end{equation*}
corresponding to $zq^4 P_{1}\bigl(zq^{5},q\bigr)$, or none of these occur
\begin{equation*}
\begin{matrix}
\cdot && \cdot && \cdot && f_5 && f_7 && f_9 && f_{11} && f_{13} & \\
& \cdot && \cdot && f_4 && f_6 && f_8 && f_{10} && f_{12} && f_{14} \\
\cdot && \cdot && f_3 && f_5 && f_7 && f_9 && f_{11} && f_{13} & \\
& \cdot && f_2 && f_4 && f_6 && f_8 && f_{10} && f_{12} && f_{14} \\
\cdot && \rs && f_3 && f_5 && f_7 && f_9 && f_{11} && f_{13} & \\
& \cdot && \rs && f_4 && f_6 && f_8 && f_{10} && f_{12} && f_{14} \\
\cdot && \cdot && \rs && f_5 && f_7 && f_9 && f_{11} && f_{13} & \\
& \cdot && \cdot && \rs && f_6 && f_8 && f_{10} && f_{12} && f_{14} \\
 \end{matrix}\quad\dots,
\end{equation*}
corresponding to $P_4(zq,q)$.
\end{itemize}

\subsection[Propositions~\ref{prop:comp1} and~\ref{prop:comp2} for [l=4]{Propositions~\ref{prop:comp1} and~\ref{prop:comp2} for $\boldsymbol{\ell=4}$}\label{subsec:sercomp}

Let us explicitly show how to obtain Propositions~\ref{prop:comp1} and~\ref{prop:comp2} from Lemma~\ref{lem:atomAG} in the case $\ell=4$. Our atomic relations are
\begin{align*}
\rel^1_{\langle v_1,v_2,v_3,v_4\rangle}:=S_{\langle v_1,v_2,v_3,v_4\rangle} - S_{\langle v_1+1,v_2,v_3,v_4\rangle} - zq^{v_1+1} S_{\langle v_1+2,v_2+2,v_3+2,v_4+2\rangle}, \\
\rel^2_{\langle v_1,v_2,v_3,v_4\rangle}:=S_{\langle v_1,v_2,v_3,v_4\rangle} - S_{\langle v_1,v_2+1,v_3,v_4\rangle} - zq^{v_2+2} S_{\langle v_1+2,v_2+4,v_3+4,v_4+4\rangle}, \\
\rel^3_{\langle v_1,v_2,v_3,v_4\rangle}:=S_{\langle v_1,v_2,v_3,v_4\rangle} - S_{\langle v_1,v_2,v_3+1,v_4\rangle} - zq^{v_3+3} S_{\langle v_1+2,v_2+4,v_3+6,v_4+6\rangle}, \\
\rel^4_{\langle v_1,v_2,v_3,v_4\rangle}:=S_{\langle v_1,v_2,v_3,v_4\rangle} - S_{\langle v_1,v_2,v_3,v_4+1\rangle} - zq^{v_4+4} S_{\langle v_1+2,v_2+4,v_3+6,v_4+8\rangle}.
\end{align*}

\begin{Proposition} The following equation is true:
\begin{equation*}S_{\langle 0,0,0,0\rangle} - S_{\langle 1,1,1,1\rangle} -
S_{\langle 1,1,1,2\rangle} + (1-zq) S_{\langle 2,2,2,2\rangle}=0.
\end{equation*}
\end{Proposition}
\begin{proof}
We begin by adding together $\rel^1_{\langle 0,0,0,0\rangle}+\rel^2_{\langle 1,0,0,0\rangle}+\rel^3_{\langle 1,1,0,0\rangle}+\rel^4_{\langle 1,1,1,0\rangle}$:
\begin{align}
&\bigl(S_{\langle 0,0,0,0\rangle} - S_{\langle 1,0,0,0\rangle} - zq^{1} S_{\langle 2,2,2,2\rangle}\bigr)
 +\bigl(S_{\langle 1,0,0,0\rangle} - S_{\langle 1,1,0,0\rangle} - zq^{2} S_{\langle 3,4,4,4\rangle}\bigr) \notag\\
& \qquad {}+\bigl(S_{\langle 1,1,0,0\rangle} - S_{\langle 1,1,1,0\rangle} - zq^{3} S_{\langle 3,5,6,6\rangle}\bigr) +\bigl(S_{\langle 1,1,1,0\rangle} - S_{\langle 1,1,1,1\rangle} - zq^{4} S_{\langle 3,5,7,8\rangle} \bigr) \notag\\
& \quad\qquad {}=S_{\langle 0,0,0,0\rangle} - S_{\langle 1,1,1,1\rangle} - zq S_{\langle 2,2,2,2\rangle} - zq^{2} S_{\langle 3,4,4,4\rangle} - zq^{3} S_{\langle 3,5,6,6\rangle}
- zq^{4} S_{\langle 3,5,7,8\rangle}.\!\!\! \label{eq:a1}
\end{align}
Meanwhile, we also add together $\rel^1_{\langle 1,2,2,2\rangle}+\rel^2_{\langle 1,1,2,2\rangle}+\rel^3_{\langle 1,1,1,2\rangle}$:
\begin{align}
&\bigl(S_{\langle 1,2,2,2\rangle} - S_{\langle 2,2,2,2\rangle} - zq^{2} S_{\langle 3,4,4,4\rangle}\bigr) +\bigl(S_{\langle 1,1,2,2\rangle} - S_{\langle 1,2,2,2\rangle} - zq^{3} S_{\langle 3,5,6,6\rangle}\bigr) \notag\\
& \qquad {}+\bigl(S_{\langle 1,1,1,2\rangle} - S_{\langle 1,1,2,2\rangle} - zq^{4} S_{\langle 3,5,7,8\rangle}\bigr) \notag\\
& \qquad\quad{} = S_{\langle 1,1,1,2\rangle}- S_{\langle 2,2,2,2\rangle} - zq^{2} S_{\langle 3,4,4,4\rangle}
- zq^{3} S_{\langle 3,5,6,6\rangle} - zq^{4} S_{\langle 3,5,7,8\rangle}. \label{eq:a2}
\end{align}
So, subtracting~\eqref{eq:a2} from~\eqref{eq:a1} produces
\begin{align*}
&S_{\langle 0,0,0,0\rangle} - S_{\langle 1,1,1,1\rangle} - zq S_{\langle 2,2,2,2\rangle} - zq^{2} S_{\langle 3,4,4,4\rangle} - zq^{3} S_{\langle 3,5,6,6\rangle}
- zq^{4} S_{\langle 3,5,7,8\rangle} \\
& \qquad\quad{} -\bigl(S_{\langle 1,1,1,2\rangle}- S_{\langle 2,2,2,2\rangle} - zq^{2} S_{\langle 3,4,4,4\rangle}
- zq^{3} S_{\langle 3,5,6,6\rangle} - zq^{4} S_{\langle 3,5,7,8\rangle}\bigr) \\
& \qquad =S_{\langle 0,0,0,0\rangle} - S_{\langle 1,1,1,1\rangle} - S_{\langle 1,1,1,2\rangle}+(1-zq)S_{\langle 2,2,2,2\rangle} =0,
\end{align*}
as desired.
\end{proof}

\begin{Proposition}
The following equations are true:
\begin{gather*}
 S_{\langle 0, 0, 0, 1\rangle}- S_{\langle 1, 1, 1, 1\rangle} - S_{\langle 1, 1, 2, 3\rangle} + (1-zq)S_{\langle 2, 2, 2, 3\rangle} =0, \\
 S_{\langle 0, 0, 1,2\rangle}- S_{\langle 1, 1, 1, 2\rangle} - S_{\langle 1, 2,3,4\rangle} + (1-zq)S_{\langle 2, 2, 3, 4\rangle} =0, \\
 S_{\langle 0, 1,2,3\rangle}- S_{\langle 1,1,2,3\rangle} - S_{\langle 2,3,4,5\rangle} + (1-zq)S_{\langle 2, 3,4,5\rangle} =0.
\end{gather*}
\end{Proposition}

\begin{proof}
We compute
\begin{align*}
&\rel^1_{\langle0,0,0,1\rangle}+\rel^2_{\langle 1,0,0,1\rangle}+\rel^3_{\langle 1,1,0,1\rangle}
-\rel^1_{\langle 1,2,2,3\rangle }-\rel^2_{\langle 1,1,2,3\rangle} \\
& \qquad {}=\bigl(S_{\langle 0, 0, 0, 1\rangle} - S_{\langle 1, 0, 0, 1\rangle} -zq S_{\langle 2, 2, 2, 3\rangle}\bigr) + \bigl(S_{\langle 1, 0, 0, 1\rangle} - S_{\langle 1, 1, 0, 1\rangle} -zq^2 S_{\langle 3, 4, 4, 5\rangle}\bigr) \\
& \qquad\quad{} +\bigl( S_{\langle 1, 1, 0, 1\rangle} - S_{\langle 1, 1, 1, 1\rangle} -zq^3 S_{\langle 3, 5, 6, 7\rangle}\bigr) -\bigl(S_{\langle 1, 2, 2, 3\rangle} - S_{\langle 2, 2, 2, 3\rangle}-zq^2 S_{\langle 3, 4, 4, 5\rangle}\bigr) \\
& \qquad \quad{} -\bigl(S_{\langle 1, 1, 2, 3\rangle} - S_{\langle 1, 2, 2, 3\rangle}-zq^3 S_{\langle 3, 5, 6, 7\rangle}\bigr) \\
&\qquad =S_{\langle 0, 0, 0, 1\rangle}- S_{\langle 1, 1, 1, 1\rangle} - S_{\langle 1, 1, 2, 3\rangle} + (1-zq)S_{\langle 2, 2, 2, 3\rangle} =0,
\\
&\rel^1_{\langle0,0,1,2\rangle}+\rel^2_{\langle 1,0,1,2\rangle}-\rel^1_{\langle 1,2,3,4\rangle}\\
&\qquad =\bigl(S_{\langle 0, 0, 1, 2\rangle} - S_{\langle 1, 0, 1, 2\rangle} -zq S_{\langle 2, 2, 3, 4\rangle}\bigr) + \bigl(S_{\langle 1, 0, 1, 2\rangle} - S_{\langle 1, 1, 1, 2\rangle} -zq^2 S_{\langle 3, 4, 5, 6\rangle}\bigr) \\
& \qquad\quad{} -\bigl(S_{\langle 1, 2, 3, 4\rangle} - S_{\langle 2, 2, 3, 4\rangle}-zq^2 S_{\langle 3, 4, 5, 6\rangle}\bigr) \\
& \qquad =S_{\langle 0, 0, 1, 2\rangle}- S_{\langle 1, 1, 1, 2\rangle} - S_{\langle 1, 2,3,4\rangle} + (1-zq)S_{\langle 2, 2, 3,4\rangle} =0,
\\
&\rel^1_{\langle0,1,2,3\rangle} =S_{\langle 0, 1, 2,3\rangle} - S_{\langle 1, 1,2,3\rangle} -zq S_{\langle 2, 3, 4, 5\rangle} \\
&\hphantom{\rel^1_{\langle0,1,2,3\rangle}} {}=S_{\langle 0, 1,2,3\rangle}- S_{\langle 1,1,2,3\rangle} - S_{\langle 2,3,4,5\rangle} + (1-zq)S_{\langle 2, 3,4,5\rangle} =0. \tag*{\qed}
\end{align*}
\renewcommand{\qed}{}
\end{proof}

\subsection[Theorem~\ref{thm:funceq1B} for l=4]{Theorem~\ref{thm:funceq1B} for $\boldsymbol{\ell=4}$}\label{subsec:combspec2}
Let us again consider $\ell=4$, resulting in $2(4)-1=7$ rows. Once again, there are five generating functions to consider:
\begin{gather*}
P^\star_4(z,q)\colon \quad
\begin{matrix}
& \cdot && \cdot && f_4 && f_6 && f_8 && f_{10} && f_{12} && f_{14} \\
 \cdot && \cdot && f_3 && f_5 && f_7 && f_9 && f_{11} && f_{13} & \\
& \cdot && f_2 && f_4 && f_6 && f_8 && f_{10} && f_{12} && f_{14} \\
\ficone && f_1 && f_3 && f_5 && f_7 && f_9 && f_{11} && f_{13} & \\
& \cdot && f_2 && f_4 && f_6 && f_8 && f_{10} && f_{12} && f_{14} \\
 \cdot && \cdot && f_3 && f_5 && f_7 && f_9 && f_{11} && f_{13} & \\
& \cdot && \cdot && f_4 && f_6 && f_8 && f_{10} && f_{12} && f_{14} \vspace{3mm}
\end{matrix}\quad\dots\,,
\\
P^\star_{3}(z,q)\colon \quad
\begin{matrix}
\cdot && \cdot && f_3 && f_5 && f_7 && f_9 && f_{11} && f_{13} & \\
& \cdot && f_2 && f_4 && f_6 && f_8 && f_{10} && f_{12} && f_{14} \\
\ficone && f_1 && f_3 && f_5 && f_7 && f_9 && f_{11} && f_{13} & \\
& \cdot && f_2 && f_4 && f_6 && f_8 && f_{10} && f_{12} && f_{14} \\
\cdot && \cdot && f_3 && f_5 && f_7 && f_9 && f_{11} && f_{13} & \\
& \cdot && \cdot && f_4 && f_6 && f_8 && f_{10} && f_{12} && f_{14} \\
\cdot && \cdot && \cdot && f_5 && f_7 && f_9 && f_{11} && f_{13} & \vspace{3mm}
\end{matrix}\quad\dots\,,
\\
P^\star_{2}(z,q)\colon \quad
\begin{matrix}
& \cdot && f_2 && f_4 && f_6 && f_8 && f_{10} && f_{12} && f_{14} \\
\ficone && f_1 && f_3 && f_5 && f_7 && f_9 && f_{11} && f_{13} & \\
& \cdot && f_2 && f_4 && f_6 && f_8 && f_{10} && f_{12} && f_{14} \\
 \cdot && \cdot && f_3 && f_5 && f_7 && f_9 && f_{11} && f_{13} & \\
& \cdot && \cdot&& f_4 && f_6 && f_8 && f_{10} && f_{12} && f_{14} \\
 \cdot && \cdot && \cdot && f_5 && f_7 && f_9 && f_{11} && f_{13} & \\
& \cdot && \cdot && \cdot && f_6 && f_8 && f_{10} && f_{12} && f_{14} \vspace{3mm}
\end{matrix}\quad\dots\,,
\\
P^\star_{1}(z,q)\colon \quad
\begin{matrix}
\ficone && f_1 && f_3 && f_5 && f_7 && f_9 && f_{11} && f_{13} & \\
& \cdot && f_2 && f_4 && f_6 && f_8 && f_{10} && f_{12} && f_{14} \\
\cdot && \cdot && f_3 && f_5 && f_7 && f_9 && f_{11} && f_{13} & \\
& \cdot && \cdot && f_4 && f_6 && f_8 && f_{10} && f_{12} && f_{14} \\
\cdot && \cdot && \cdot && f_5 && f_7 && f_9 && f_{11} && f_{13} & \\
& \cdot && \cdot && \cdot && f_6 && f_8 && f_{10} && f_{12} && f_{14} \\
\cdot && \cdot && \cdot && \cdot && f_7 && f_9 && f_{11} && f_{13} & \vspace{3mm}
\end{matrix}\quad\dots\,,
\\
P^\star_{0}(z,q)\colon \quad
\begin{matrix}
& \ficone && f_2 && f_4 && f_6 && f_8 && f_{10} && f_{12} && f_{14} \\
\cdot && \cdot && f_3 && f_5 && f_7 && f_9 && f_{11} && f_{13} & \\
& \cdot && \cdot && f_4 && f_6 && f_8 && f_{10} && f_{12} && f_{14} \\
\cdot && \cdot && \cdot && f_5 && f_7 && f_9 && f_{11} && f_{13} & \\
& \cdot && \cdot && \cdot && f_6 && f_8 && f_{10} && f_{12} && f_{14} \\
\cdot && \cdot && \cdot && \cdot && f_7 && f_9 && f_{11} && f_{13} & \\
& \cdot && \cdot && \cdot && \cdot && f_8 && f_{10} && f_{12} && f_{14}
\end{matrix}\quad\dots\,.
\end{gather*}
\begin{itemize}\itemsep=0pt
\item $P^\star_{0}(z,q) = P^\star_{1}(zq,q)$: This is clear by inspection.

\item $P^\star_{1}(z,q) =
zq^1 P^\star_{1}\bigl(zq^{2},q\bigr)
+ P^\star_{2}(zq,q)$:
Either there is a 1 in the first row
\begin{equation*}
\begin{matrix}
\cdot && \color{red}{1} && f_3 && f_5 && f_7 && f_9 && f_{11} && f_{13} & \\
& \cdot && \rs&& f_4 && f_6 && f_8 && f_{10} && f_{12} && f_{14} \\
\cdot && \cdot && \rs && f_5 && f_7 && f_9 && f_{11} && f_{13} & \\
& \cdot && \cdot && \rs && f_6 && f_8 && f_{10} && f_{12} && f_{14} \\
\cdot && \cdot && \cdot && \rs && f_7 && f_9 && f_{11} && f_{13} & \\
& \cdot && \cdot && \cdot && \rs && f_8 && f_{10} && f_{12} && f_{14} \\
\cdot && \cdot && \cdot && \cdot && \rs && f_9 && f_{11} && f_{13} &
\end{matrix}\quad\dots\,,
\end{equation*}
which corresponds to $zq^1 P^\star_{1}\bigl(zq^{2},q\bigr)$, or there is no 1
\begin{equation*}
\begin{matrix}
\cdot && \rs && f_3 && f_5 && f_7 && f_9 && f_{11} && f_{13} & \\
& \cdot && f_2 && f_4 && f_6 && f_8 && f_{10} && f_{12} && f_{14} \\
\cdot && \cdot && f_3 && f_5 && f_7 && f_9 && f_{11} && f_{13} & \\
& \cdot && \cdot && f_4 && f_6 && f_8 && f_{10} && f_{12} && f_{14} \\
\cdot && \cdot && \cdot && f_5 && f_7 && f_9 && f_{11} && f_{13} & \\
& \cdot && \cdot && \cdot && f_6 && f_8 && f_{10} && f_{12} && f_{14} \\
\cdot && \cdot && \cdot && \cdot && f_7 && f_9 && f_{11} && f_{13} &
\end{matrix}\quad\dots\,,
\end{equation*}
which corresponds to $P^\star_{2}(zq,q)$.

\item $P^\star_{2}(z,q) =
zq^1 P^\star_{2}\bigl(zq^{2},q\bigr) +
zq^2 P^\star_{1}\bigl(zq^{3},q\bigr) +
P^\star_{3}(zq,q)$: Either there is a 1
\begin{equation*}
\begin{matrix}
& \cdot && \rs && f_4 && f_6 && f_8 && f_{10} && f_{12} && f_{14} \\
\cdot && \color{red}{1} && f_3 && f_5 && f_7 && f_9 && f_{11} && f_{13} & \\
& \cdot && \rs && f_4 && f_6 && f_8 && f_{10} && f_{12} && f_{14} \\
 \cdot && \cdot && \rs && f_5 && f_7 && f_9 && f_{11} && f_{13} & \\
& \cdot && \cdot&& \rs && f_6 && f_8 && f_{10} && f_{12} && f_{14} \\
 \cdot && \cdot && \cdot && \rs && f_7 && f_9 && f_{11} && f_{13} & \\
& \cdot && \cdot && \cdot && \rs && f_8 && f_{10} && f_{12} && f_{14}
\end{matrix}\quad\dots\,,
\end{equation*}
corresponding to $zq^1 P^\star_{2}\bigl(zq^{2},q\bigr)$, or there is a 2 in the first row
\begin{equation*}
\begin{matrix}
& \cdot && \color{red}{2} && f_4 && f_6 && f_8 && f_{10} && f_{12} && f_{14} \\
\cdot && \rs && \rs && f_5 && f_7 && f_9 && f_{11} && f_{13} & \\
& \cdot && \rs && \rs && f_6 && f_8 && f_{10} && f_{12} && f_{14} \\
 \cdot && \cdot && \rs && \rs && f_7 && f_9 && f_{11} && f_{13} & \\
& \cdot && \cdot&& \rs && \rs && f_8 && f_{10} && f_{12} && f_{14} \\
 \cdot && \cdot && \cdot && \rs && \rs && f_9 && f_{11} && f_{13} & \\
& \cdot && \cdot && \cdot && \rs && \rs && f_{10} && f_{12} && f_{14}
\end{matrix}\quad\dots\,,
\end{equation*}
corresponding to $zq^2 P^\star_{1}\bigl(zq^{3},q\bigr) $, or neither of these occur
\begin{equation*}
\begin{matrix}
& \cdot && \rs && f_4 && f_6 && f_8 && f_{10} && f_{12} && f_{14} \\
\cdot && \rs && f_3 && f_5 && f_7 && f_9 && f_{11} && f_{13} & \\
& \cdot && f_2 && f_4 && f_6 && f_8 && f_{10} && f_{12} && f_{14} \\
 \cdot && \cdot && f_3 && f_5 && f_7 && f_9 && f_{11} && f_{13} & \\
& \cdot && \cdot&& f_4 && f_6 && f_8 && f_{10} && f_{12} && f_{14} \\
 \cdot && \cdot && \cdot && f_5 && f_7 && f_9 && f_{11} && f_{13} & \\
& \cdot && \cdot && \cdot && f_6 && f_8 && f_{10} && f_{12} && f_{14}
\end{matrix}\quad\dots\,,
\end{equation*}
corresponding to $P^\star_{3}(zq,q)$.

\item $P^\star_{3}(z,q) =
zq^1 P^\star_{3}\bigl(zq^{2},q\bigr) +
zq^2 P^\star_{2}\bigl(zq^{3},q\bigr) +
zq^3 P^\star_{1}\bigl(zq^{4},q\bigr) +
P^\star_{4}(zq,q)$:
Either there is a 1
\begin{equation*}
\begin{matrix}
\cdot && \cdot && \rs && f_5 && f_7 && f_9 && f_{11} && f_{13} & \\
& \cdot && \rs && f_4 && f_6 && f_8 && f_{10} && f_{12} && f_{14} \\
\cdot && \color{red}{1} && f_3 && f_5 && f_7 && f_9 && f_{11} && f_{13} & \\
& \cdot && \rs && f_4 && f_6 && f_8 && f_{10} && f_{12} && f_{14} \\
\cdot && \cdot && \rs && f_5 && f_7 && f_9 && f_{11} && f_{13} & \\
& \cdot && \cdot && \rs && f_6 && f_8 && f_{10} && f_{12} && f_{14} \\
\cdot && \cdot && \cdot && \rs && f_7 && f_9 && f_{11} && f_{13} &
\end{matrix}\quad\dots\,,
\end{equation*}
corresponding to $zq^1 P^\star_{3}\bigl(zq^{2},q\bigr)$, or there is a 2 in the second row
\begin{equation*}
\begin{matrix}
\cdot && \cdot && \rs && f_5 && f_7 && f_9 && f_{11} && f_{13} & \\
& \cdot && \color{red}{2} && f_4 && f_6 && f_8 && f_{10} && f_{12} && f_{14} \\
\cdot && \rs && \rs && f_5 && f_7 && f_9 && f_{11} && f_{13} & \\
& \cdot && \rs && \rs && f_6 && f_8 && f_{10} && f_{12} && f_{14} \\
\cdot && \cdot && \rs && \rs && f_7 && f_9 && f_{11} && f_{13} & \\
& \cdot && \cdot && \rs && \rs && f_8 && f_{10} && f_{12} && f_{14} \\
\cdot && \cdot && \cdot && \rs && \rs && f_9 && f_{11} && f_{13} &
\end{matrix}\quad\dots\,,
\end{equation*}
corresponding to $zq^2 P^\star_{2}\bigl(zq^{3},q\bigr)$, or there is a 3 in the first row
\begin{equation*}
\begin{matrix}
\cdot && \cdot && \color{red}{3} && f_5 && f_7 && f_9 && f_{11} && f_{13} & \\
& \cdot && \rs && \rs && f_6 && f_8 && f_{10} && f_{12} && f_{14} \\
\cdot && \rs && \rs && \rs && f_7 && f_9 && f_{11} && f_{13} & \\
& \cdot && \rs && \rs && \rs && f_8 && f_{10} && f_{12} && f_{14} \\
\cdot && \cdot && \rs && \rs && \rs && f_9 && f_{11} && f_{13} & \\
& \cdot && \cdot && \rs && \rs && \rs && f_{10} && f_{12} && f_{14} \\
\cdot && \cdot && \cdot && \rs && \rs && \rs && f_{11} && f_{13} &
\end{matrix}\quad\dots\,,
\end{equation*}
corresponding to $zq^3 P^\star_{1}\bigl(zq^{4},q\bigr)$, or none of these occur
\begin{equation*}
\begin{matrix}
\cdot && \cdot && \rs && f_5 && f_7 && f_9 && f_{11} && f_{13} & \\
& \cdot && \rs && f_4 && f_6 && f_8 && f_{10} && f_{12} && f_{14} \\
\cdot && \rs && f_3 && f_5 && f_7 && f_9 && f_{11} && f_{13} & \\
& \cdot && f_2 && f_4 && f_6 && f_8 && f_{10} && f_{12} && f_{14} \\
\cdot && \cdot && f_3 && f_5 && f_7 && f_9 && f_{11} && f_{13} & \\
& \cdot && \cdot && f_4 && f_6 && f_8 && f_{10} && f_{12} && f_{14} \\
\cdot && \cdot && \cdot && f_5 && f_7 && f_9 && f_{11} && f_{13} &
\end{matrix}\quad\dots\,,
\end{equation*}
corresponding to $P^\star_{4}(zq,q)$.

\item $P^\star_4(z,q) =
zq^1 P^\star_{4}\bigl(zq^{2},q\bigr) +
zq^2 P^\star_{3}\bigl(zq^{3},q\bigr) +
zq^3 P^\star_{2}\bigl(zq^{4},q\bigr) +
zq^4 P^\star_{1}\bigl(zq^{5},q\bigr)+
P^\star_{3}(zq,q)$:
Either there is a~1
\begin{equation*}
\begin{matrix}
& \cdot && \cdot && \rs && f_6 && f_8 && f_{10} && f_{12} && f_{14} \\
 \cdot && \cdot && \rs && f_5 && f_7 && f_9 && f_{11} && f_{13} & \\
& \cdot && \rs && f_4 && f_6 && f_8 && f_{10} && f_{12} && f_{14} \\
\cdot && \color{red}{1} && f_3 && f_5 && f_7 && f_9 && f_{11} && f_{13} & \\
& \cdot && \rs && f_4 && f_6 && f_8 && f_{10} && f_{12} && f_{14} \\
 \cdot && \cdot && \rs && f_5 && f_7 && f_9 && f_{11} && f_{13} & \\
& \cdot && \cdot && \rs && f_6 && f_8 && f_{10} && f_{12} && f_{14}
\end{matrix}\quad\dots\,,
\end{equation*}
corresponding to $zq^1 P^\star_{4}\bigl(zq^{2},q\bigr)$, or there is a 2 in the third row
\begin{equation*}
\begin{matrix}
& \cdot && \cdot && \rs && f_6 && f_8 && f_{10} && f_{12} && f_{14} \\
 \cdot && \cdot && \rs && f_5 && f_7 && f_9 && f_{11} && f_{13} & \\
& \cdot && \color{red}{2} && f_4 && f_6 && f_8 && f_{10} && f_{12} && f_{14} \\
\cdot && \rs&& \rs && f_5 && f_7 && f_9 && f_{11} && f_{13} & \\
& \cdot && \rs && \rs && f_6 && f_8 && f_{10} && f_{12} && f_{14} \\
 \cdot && \cdot && \rs && \rs && f_7 && f_9 && f_{11} && f_{13} & \\
& \cdot && \cdot && \rs&& \rs && f_8 && f_{10} && f_{12} && f_{14}
\end{matrix}\quad\dots\,,
\end{equation*}
corresponding to $zq^2 P^\star_{3}\bigl(zq^{3},q\bigr)$, or there is a 3 in the second row
\begin{equation*}
\begin{matrix}
& \cdot && \cdot && \rs && f_6 && f_8 && f_{10} && f_{12} && f_{14} \\
 \cdot && \cdot && \color{red}{3} && f_5 && f_7 && f_9 && f_{11} && f_{13} & \\
& \cdot && \rs && \rs && f_6 && f_8 && f_{10} && f_{12} && f_{14} \\
\cdot && \rs && \rs && \rs && f_7 && f_9 && f_{11} && f_{13} & \\
& \cdot && \rs && \rs && \rs && f_8 && f_{10} && f_{12} && f_{14} \\
 \cdot && \cdot && \rs && \rs && f_7 && f_9 && f_{11} && f_{13} & \\
& \cdot && \cdot && \rs && \rs && f_8 && f_{10} && f_{12} && f_{14}
\end{matrix}\quad\dots\,,
\end{equation*}
corresponding to $zq^3 P^\star_{2}\bigl(zq^{4},q\bigr)$, or there is a 4 in the first row
\begin{equation*}
\begin{matrix}
& \cdot && \cdot && \color{red}{4} && f_6 && f_8 && f_{10} && f_{12} && f_{14} \\
 \cdot && \cdot && \rs && \rs && f_7 && f_9 && f_{11} && f_{13} & \\
& \cdot && \rs && \rs && \rs && f_8 && f_{10} && f_{12} && f_{14} \\
\cdot && \rs && \rs && \rs && \rs && f_9 && f_{11} && f_{13} & \\
& \cdot && \rs && \rs && \rs && \rs && f_{10} && f_{12} && f_{14} \\
 \cdot && \cdot && \rs && \rs && \rs && \rs && f_{11} && f_{13} & \\
& \cdot && \cdot && \rs && \rs && \rs && \rs && f_{12} && f_{14}
\end{matrix}\quad\dots\,,
\end{equation*}
corresponding to $zq^4 P^\star_{1}\bigl(zq^{5},q\bigr)$, or none of these occur
\begin{equation*}
\begin{matrix}
& \cdot && \cdot && \rs && f_6 && f_8 && f_{10} && f_{12} && f_{14} \\
 \cdot && \cdot && \rs && f_5 && f_7 && f_9 && f_{11} && f_{13} & \\
& \cdot && \rs && f_4 && f_6 && f_8 && f_{10} && f_{12} && f_{14} \\
\cdot && \rs && f_3 && f_5 && f_7 && f_9 && f_{11} && f_{13} & \\
& \cdot && f_2 && f_4 && f_6 && f_8 && f_{10} && f_{12} && f_{14} \\
 \cdot && \cdot && f_3 && f_5 && f_7 && f_9 && f_{11} && f_{13} & \\
& \cdot && \cdot && f_4 && f_6 && f_8 && f_{10} && f_{12} && f_{14}
\end{matrix}\quad\dots\,,
\end{equation*}
corresponding to $P^\star_{3}(zq,q)$.
\end{itemize}

\subsection*{Acknowledgements}

The author would like to express his gratitude to Shashank Kanade for many fruitful conversations and extensive feedback, to Stefano Capparelli, Arne Meurman, and Mirko Primc for corrections and encouraging comments on drafts of this paper, to Ole Warnaar for his key insights that the multisums of Theorems \ref{thm:mainresult1} and~\ref{thm:mainresult2} are (up to a factor of $(zq;q)_\infty$) the same as the generating functions for two-rowed cylindric partitions and that the functional equations in Sections~\ref{sec:comp} and~\ref{sec:comp5} are the same as the (normalized) Corteel--Welsh equations, to Shunsuke Tsuchioka for (as previously noted) kindly pointing out that the $k=1$ cases of these identities are initially due (in disguised form) to Jing, Misra, and Savage~\cite{JMS}. The feedback provided by the anonymous referees was also extremely valuable, especially as it resulted in improvements in the proofs in Sections~\ref{sec:partitions2},~\ref{sec:partitions4}, and~\ref{sec:6}.


\begin{thebibliography}{99}
\footnotesize\itemsep=0pt

\bibitem{ADJM}
Afsharijoo P., Dousse J., Jouhet F., Mourtada H., New companions to the
 {A}ndrews--{G}ordon identities motivated by commutative algebra,
 \href{https://doi.org/10.1016/j.aim.2023.108946}{\textit{Adv. Math.}}
 \textbf{417} (2023), 108946, 40~pages,
 \href{http://arxiv.org/abs/2104.09422}{arXiv:2104.09422}.

\bibitem{And66}
Andrews G.E., An analytic proof of the {R}ogers--{R}amanujan--{G}ordon
 identities, \href{https://doi.org/10.2307/2373082}{\textit{Amer.~J. Math.}}
 \textbf{88} (1966), 844--846.

\bibitem{And_GGA1}
Andrews G.E., A generalization of the {G}\"{o}llnitz--{G}ordon partition
 theorems, \href{https://doi.org/10.2307/2035143}{\textit{Proc. Amer. Math.
 Soc.}} \textbf{18} (1967), 945--952.

\bibitem{AndPNAS}
Andrews G.E., An analytic generalization of the {R}ogers--{R}amanujan
 identities for odd moduli,
 \href{https://doi.org/10.1073/pnas.71.10.4082}{\textit{Proc. Nat. Acad. Sci.
 USA}} \textbf{71} (1974), 4082--4085.

\bibitem{And_GGA2}
Andrews G.E., Problems and prospects for basic hypergeometric functions, in
 Theory and {A}pplication of {S}pecial {F}unctions,
 \href{https://doi.org/10.1016/B978-0-12-064850-4.50008-2}{Academic Press}, New
 York, 1975, 191--224.

\bibitem{AndDurf}
Andrews G.E., Partitions and {D}urfee dissection,
 \href{https://doi.org/10.2307/2373804}{\textit{Amer.~J. Math.}} \textbf{101}
 (1979), 735--742.

\bibitem{And-book}
Andrews G.E., The theory of partitions, \textit{Cambridge Math. Lib.},
 \href{https://doi.org/10.1017/CBO9780511608650}{Cambridge University Press},
 Cambridge, 1998.

\bibitem{AndParity}
Andrews G.E., Parity in partition identities,
 \href{https://doi.org/10.1007/s11139-008-9150-0}{\textit{Ramanujan~J.}}
 \textbf{23} (2010), 45--90.

\bibitem{BKRS}
Baker K., Kanade S., Russell M.C., Sadowski C., Principal subspaces of basic
 modules for twisted affine {L}ie algebras, {$q$}-series multisums, and
 {N}andi's identities, \href{https://doi.org/10.5802/alco.311}{\textit{Algebr.
 Comb.}} \textbf{6} (2023), 1533--1556,
 \href{http://arxiv.org/abs/2208.14581}{arXiv:2208.14581}.

\bibitem{Bor}
Borodin A., Periodic {S}chur process and cylindric partitions,
 \href{https://doi.org/10.1215/S0012-7094-07-14031-6}{\textit{Duke Math.~J.}}
 \textbf{140} (2007), 391--468,
 \href{http://arxiv.org/abs/math.CO/0601019}{arXiv:math.CO/0601019}.

\bibitem{Bress0}
Bressoud D.M., A generalization of the {R}ogers--{R}amanujan identities for all
 moduli,
 \href{https://doi.org/10.1016/0097-3165(79)90008-6}{\textit{J.~Combin. Theory
 Ser.~A}} \textbf{27} (1979), 64--68.

\bibitem{Bress1}
Bressoud D.M., Analytic and combinatorial generalizations of the
 {R}ogers--{R}amanujan identities,
 \href{https://doi.org/10.1090/memo/0227}{\textit{Mem. Amer. Math. Soc.}}
 \textbf{24} (1980), 54~pages.

\bibitem{Bress2}
Bressoud D.M., An analytic generalization of the {R}ogers--{R}amanujan
 identities with interpretation,
 \href{https://doi.org/10.1093/qmath/31.4.385}{\textit{Q.~J. Math.}}
 \textbf{31} (1980), 385--399.

\bibitem{BZ}
Bressoud D.M., Zeilberger D., A short {R}ogers--{R}amanujan bijection,
 \href{https://doi.org/10.1016/0012-365X(82)90298-9}{\textit{Discrete Math.}}
 \textbf{38} (1982), 313--315.

\bibitem{Capp}
Capparelli S., A construction of the level {$3$} modules for the affine {L}ie
 algebra \smash{$A^{(2)}_2$} and a new combinatorial identity of the
 {R}ogers--{R}amanujan type,
 \href{https://doi.org/10.1090/S0002-9947-96-01535-8}{\textit{Trans. Amer.
 Math. Soc.}} \textbf{348} (1996), 481--501.

\bibitem{CMPP}
Capparelli S., Meurman A., Primc A., Primc M., New partition identities from
 \smash{$C^{(1)}_\ell$}-modules,
 \href{https://doi.org/10.3336/gm.57.2.01}{\textit{Glas. Mat. Ser.~III}}
 \textbf{57(77)} (2022), 161--184,
 \href{http://arxiv.org/abs/2106.06262}{arXiv:2106.06262}.

\bibitem{Cort}
Corteel S., Rogers--{R}amanujan identities and the
 {R}obinson--{S}chensted--{K}nuth correspondence,
 \href{https://doi.org/10.1090/proc/13373}{\textit{Proc. Amer. Math. Soc.}}
 \textbf{145} (2017), 2011--2022.

\bibitem{CW}
Corteel S., Welsh T., The {$A_2$} {R}ogers--{R}amanujan identities revisited,
 \href{https://doi.org/10.1007/s00026-019-00446-7}{\textit{Ann. Comb.}}
 \textbf{23} (2019), 683--694,
 \href{http://arxiv.org/abs/1905.08343}{arXiv:1905.08343}.

\bibitem{CKLMQRS}
Coulson B., Kanade S., Lepowsky J., McRae R., Qi F., Russell M.C., Sadowski C.,
 A motivated proof of the {G}\"{o}llnitz--{G}ordon--{A}ndrews identities,
 \href{https://doi.org/10.1007/s11139-015-9722-8}{\textit{Ramanujan~J.}}
 \textbf{42} (2017), 97--129,
 \href{http://arxiv.org/abs/1411.2044}{arXiv:1411.2044}.

\bibitem{DouKon}
Dousse J., Konan I., Characters of level $1$ standard modules of \smash{${C}_n^{(1)}$}
 as generating functions for generalised partitions,
 \href{http://arxiv.org/abs/2212.12728}{arXiv:2212.12728}.

\bibitem{FodWel}
Foda O., Welsh T.A., Cylindric partitions, {$\mathcal{W}_r$} characters and the
 {A}ndrews--{G}ordon--{B}ressoud identities,
 \href{https://doi.org/10.1088/1751-8113/49/16/164004}{\textit{J.~Phys.~A}}
 \textbf{49} (2016), 164004, 37~pages,
 \href{http://arxiv.org/abs/1510.02213}{arXiv:1510.02213}.

\bibitem{GarMil}
Garsia A.M., Milne S.C., A {R}ogers--{R}amanujan bijection,
 \href{https://doi.org/10.1016/0097-3165(81)90062-5}{\textit{J.~Combin. Theory
 Ser.~A}} \textbf{31} (1981), 289--339.

\bibitem{Gar}
Garvan F., {A} {$q$}-product tutorial for a
 {$q$}-series {MAPLE} package, \textit{S\'{e}m. Lothar. Combin.} \textbf{42}
 (1999), B42d, 27~pages, \url{https://www.mat.univie.ac.at/~slc/wpapers/s42garvan.html},
 \href{http://arxiv.org/abs/math.CO/9812092}{arXiv:math.CO/9812092}.

\bibitem{GesKra}
Gessel I.M., Krattenthaler C., Cylindric partitions,
 \href{https://doi.org/10.1090/S0002-9947-97-01791-1}{\textit{Trans. Amer.
 Math. Soc.}} \textbf{349} (1997), 429--479.

\bibitem{Gord}
Gordon B., A combinatorial generalization of the {R}ogers--{R}amanujan
 identities, \href{https://doi.org/10.2307/2372962}{\textit{Amer.~J. Math.}}
 \textbf{83} (1961), 393--399.

\bibitem{JMS}
Jing N., Misra K.C., Savage C.D., On multi-color partitions and the generalized
 {R}ogers--{R}amanujan identities,
 \href{https://doi.org/10.1142/S0219199701000482}{\textit{Commun. Contemp.
 Math.}} \textbf{3} (2001), 533--548,
 \href{http://arxiv.org/abs/math.CO/9907183}{arXiv:math.CO/9907183}.

\bibitem{KLRS}
Kanade S., Lepowsky J., Russell M.C., Sills A.V., Ghost series and a motivated
 proof of the {A}ndrews--{B}ressoud identities,
 \href{https://doi.org/10.1016/j.jcta.2016.07.004}{\textit{J.~Combin. Theory
 Ser.~A}} \textbf{146} (2017), 33--62,
 \href{http://arxiv.org/abs/1411.2048}{arXiv:1411.2048}.

\bibitem{KanRus-idf}
Kanade S., Russell M.C., {\tt {I}dentity{F}inder} and some new identities of
 {R}ogers--{R}amanujan type,
 \href{https://doi.org/10.1080/10586458.2015.1015186}{\textit{Exp. Math.}}
 \textbf{24} (2015), 419--423,
 \href{http://arxiv.org/abs/1411.5346}{arXiv:1411.5346}.

\bibitem{KanRus-stair}
Kanade S., Russell M.C., Staircases to analytic sum-sides for many new integer
 partition identities of {R}ogers--{R}amanujan type,
 \href{https://doi.org/10.37236/7847}{\textit{Electron.~J. Combin.}}
 \textbf{26} (2019), 1.6, 33~pages,
 \href{http://arxiv.org/abs/1803.02515}{arXiv:1803.02515}.

\bibitem{KanRus-cylindric}
Kanade S., Russell M.C., Completing the {$\mathrm{A}_2$}
 {A}ndrews--{S}chilling--{W}arnaar identities,
 \href{https://doi.org/10.1093/imrn/rnac217}{\textit{Int. Math. Res. Not.}}
 \textbf{2023} (2023), 17100--17155,
 \href{http://arxiv.org/abs/2203.05690}{arXiv:2203.05690}.

\bibitem{KRTW}
Kanade S., Russell M.C., Tsuchioka S., Warnaar S.O., Remarks on the conjectures
 of {C}apparelli, {M}eurman, {P}rimc and {P}rimc, \textit{Selecta
 Math.~(N.S.)}, {t}o appear,
 \href{http://arxiv.org/abs/2404.03851}{arXiv:2404.03851}.

\bibitem{KimYee}
Kim S., Yee A.J., The {R}ogers--{R}amanujan--{G}ordon identities, the
 generalized {G}\"{o}llnitz--{G}ordon identities, and parity questions,
 \href{https://doi.org/10.1016/j.jcta.2013.02.005}{\textit{J.~Combin. Theory
 Ser.~A}} \textbf{120} (2013), 1038--1056.

\bibitem{KurThesis}
Kur\c{s}ung\"{o}z K., Parity considerations in {A}ndrews--{G}ordon identities,
 and the $k$-marked {D}urfee symbols, Ph.D.~Thesis, {T}he Pennsylvania State
 University, 2009, \url{https://www.proquest.com/docview/304983307}.

\bibitem{KurParity}
Kur\c{s}ung\"{o}z K., Parity considerations in {A}ndrews--{G}ordon identities,
 \href{https://doi.org/10.1016/j.ejc.2009.06.002}{\textit{European~J.
 Combin.}} \textbf{31} (2010), 976--1000.

\bibitem{LepWil-struI}
Lepowsky J., Wilson R.L., The structure of standard modules. {I}. {U}niversal
 algebras and the {R}ogers--{R}amanujan identities,
 \href{https://doi.org/10.1007/BF01388447}{\textit{Invent. Math.}} \textbf{77}
 (1984), 199--290.

\bibitem{MSZ}
Mc~Laughlin J., Sills A.V., Zimmer P., Rogers--{R}amanujan--{S}later type
 identities, \href{https://doi.org/10.37236/36}{\textit{Electron.~J. Combin.}}
 \textbf{DS15} (2008), 59~pages,
 \href{http://arxiv.org/abs/1901.00946}{arXiv:1901.00946}.

\bibitem{MP87}
Meurman A., Primc M., Annihilating ideals of standard modules of
 {$\mathrm{sl}(2, \mathbf{C})^{\sim}$} and combinatorial identities,
 \href{https://doi.org/10.1016/0001-8708(87)90008-9}{\textit{Adv. Math.}}
 \textbf{64} (1987), 177--240.

\bibitem{MP}
Meurman A., Primc M., Annihilating fields of standard modules of {${\mathfrak
 s}{\mathfrak l}(2,\mathbf{C})^\sim$} and combinatorial identities,
 \href{https://doi.org/10.1090/memo/0652}{\textit{Mem. Amer. Math. Soc.}}
 \textbf{137} (1999), viii+89~pages,
 \href{http://arxiv.org/abs/math.QA/9806105}{arXiv:math.QA/9806105}.

\bibitem{Nan-thesis}
Nandi D., Partition identities arising from the standard \smash{$A_2^{(2)}$}-modules of
 level~4, Ph.D.~Thesis, {T}he State University of New Jersey, 2014,
 \url{https://www.proquest.com/docview/1655001196}.

\bibitem{PakSurvey}
Pak I., Partition bijections, a survey,
 \href{https://doi.org/10.1007/s11139-006-9576-1}{\textit{Ramanujan~J.}}
 \textbf{12} (2006), 5--75.

\bibitem{primc2023}
Primc M., New partition identities for odd {$W$} odd,
 \href{https://doi.org/10.21857/9e31lhzl8m}{\textit{Rad Hrvat. Akad. Znan.
 Umjet. Mat. Znan.}} \textbf{28(558)} (2024), 49--56,
 \href{http://arxiv.org/abs/2301.12484}{arXiv:2301.12484}.

\bibitem{PT}
Primc M., Trup\v{c}evi\'{c} G., Linear independence for \smash{$C_\ell^{(1)}$} by
 using \smash{$C_{2\ell}^{(1)}$},
 \href{https://doi.org/10.1016/j.jalgebra.2024.08.003}{\textit{J.~Algebra}}
 \textbf{661} (2025), 341--356,
 \href{http://arxiv.org/abs/2403.06881}{arXiv:2403.06881}.

\bibitem{PS2}
Primc M., \v{S}iki\'{c} T., Combinatorial bases of basic modules for affine
 {L}ie algebras \smash{$C_n^{(1)}$},
 \href{https://doi.org/10.1063/1.4962392}{\textit{J.~Math. Phys.}} \textbf{57}
 (2016), 091701, 19~pages,
 \href{http://arxiv.org/abs/1603.04399}{arXiv:1603.04399}.

\bibitem{PS1}
Primc M., \v{S}iki\'{c} T., Leading terms of relations for standard modules of
 the affine {L}ie algebras \smash{$C_n^{(1)}$},
 \href{https://doi.org/10.1007/s11139-018-0052-5}{\textit{Ramanujan~J.}}
 \textbf{48} (2019), 509--543,
 \href{http://arxiv.org/abs/1506.05026}{arXiv:1506.05026}.

\bibitem{Sil-book}
Sills A.V., An invitation to the {R}ogers--{R}amanujan identities,
 \href{https://doi.org/10.1201/9781315151922}{CRC Press}, Boca Raton, FL,
 2018.

\bibitem{Stant}
Stanton D., Binomial {A}ndrews--{G}ordon--{B}ressoud identities, in Frontiers
 in {O}rthogonal {P}olynomials and {$q$}-{S}eries, \textit{Contemp. Math.
 Appl. Monogr. Expo. Lect. Notes}, Vol.~1,
 \href{https://doi.org/10.1142/9789813228887_0002}{World Scientific
 Publishing}, Hackensack, NJ, 2018, 7--19,
 \href{http://arxiv.org/abs/1608.01294}{arXiv:1608.01294}.

\bibitem{TakTsu-nandi}
Takigiku M., Tsuchioka S., A proof of conjectured partition identities of
 {N}andi, \href{https://doi.org/10.1353/ajm.2024.a923238}{\textit{Amer.~J.
 Math.}} \textbf{146} (2024), 405--433,
 \href{http://arxiv.org/abs/1910.12461}{arXiv:1910.12461}.

\bibitem{Tru}
Trup\v{c}evi\'{c} G., Bases of standard modules for affine {L}ie algebras of
 type \smash{$C^{(1)}_{\ell}$},
 \href{https://doi.org/10.1080/00927872.2018.1424874}{\textit{Comm. Algebra}}
 \textbf{46} (2018), 3663--3673,
 \href{http://arxiv.org/abs/1706.10089}{arXiv:1706.10089}.

\bibitem{WarnCyl}
Warnaar S.O., The {${\rm A}_2$} {A}ndrews--{G}ordon identities and cylindric
 partitions, \href{https://doi.org/10.1090/btran/147}{\textit{Trans. Amer.
 Math. Soc. Ser.~B}} \textbf{10} (2023), 715--765,
 \href{http://arxiv.org/abs/2111.07550}{arXiv:2111.07550}.

\bibitem{Wiecz}
Wieczorek M., Andrews--{B}ressoud series and {W}ronskians,
 \href{http://arxiv.org/abs/2002.07846}{arXiv:2002.07846}.

\end{thebibliography}

\pdfbookmark[1]{References}{ref}
\LastPageEnding

\end{document}